\documentclass[12pt]{amsart}
\usepackage{epsfig, graphics, psfrag, color, hyperref}

\title{Algorithmically detecting the bridge number
of hyperbolic knots}
\author{Alexander Coward}

\theoremstyle{plain}
\newtheorem{theorem}{Theorem}[section]
\newtheorem{corollary}[theorem]{Corollary}
\newtheorem{cor}[theorem]{Corollary}
\newtheorem{lemma}[theorem]{Lemma}
\newtheorem{prop}[theorem]{Proposition}
\newtheorem{proposition}[theorem]{Proposition}

\newtheorem*{claimx}{Claim}

\newtheorem*{namedtheorem}{\theoremname}
\newcommand{\theoremname}{Theorem \ref{untelprop}}

\newtheorem*{namedtheoremb}{\theoremnameb}
\newcommand{\theoremnameb}{Lemma \ref{pushlemma}}

\newtheorem*{namedtheoremc}{\theoremnamec}
\newcommand{\theoremnamec}{Theorem \ref{paratori}}

\newtheorem*{namedtheoremd}{\theoremnamed}
\newcommand{\theoremnamed}{Theorem \ref{21lack}}

\newtheorem*{namedtheoreme}{\theoremnamee}
\newcommand{\theoremnamee}{Theorem \ref{nktheorem}}

\newtheorem*{namedtheoremf}{\theoremnamef}
\newcommand{\theoremnamef}{Theorem \ref{findthem}}

\newtheorem*{namedtheoremg}{\theoremnameg}
\newcommand{\theoremnameg}{Theorem \ref{newthm}}

\theoremstyle{definition}
\newtheorem{define}[theorem]{Definition}

\newtheorem*{remark}{Remark}

\psfrag{F0}[][]{$F_0$} \psfrag{F1}[][]{$F_1$} \psfrag{N1}[][]{$N_1$}
\psfrag{f}[][]{$f$} \psfrag{hd}[][]{$h'$} \psfrag{I}[][]{$I$}
\psfrag{Fri}[][]{$f_i$} \psfrag{D}[][]{$D$} \psfrag{N_1}[][]{$N_1$}
\psfrag{N_2}[][]{$N_2$} \psfrag{N_3}[][]{$N_3$} \psfrag{m}[][]{$m$}
\psfrag{T'}[][]{$T'$} \psfrag{bM}[][]{$\partial M$}
\psfrag{fi}[][]{$f_i$} \psfrag{a}[cl][bl]{$a$}
\psfrag{b}[cl][bl]{$b$} \psfrag{<=a*2}[cl][bl]{$\leq 2a$}
\psfrag{<=b*2}[cl][bl]{$\leq 2b$} \psfrag{<=ab}[cl][bl]{$\leq
2(a+b)$}

\psfrag{l1}[][bl]{$I_1$} \psfrag{l2}[][bl]{$I_2$}
\psfrag{l3}[][bl]{$E_1$} \psfrag{l4}[][bl]{$E_p$}
\psfrag{l5}[][bl]{$E_1'$} \psfrag{l6}[][bl]{$E_p'$}
\psfrag{l7}[][bl][2][90]{$-$} \psfrag{l8}[][bl]{$E_{p+1}$}
\psfrag{l9}[][bl]{$E_{q}$} \psfrag{l10}[][bl]{$E_{p+2}$}
\psfrag{l11}[][bl]{$E_{q-1}$} \psfrag{l12}[][bl]{$\ldots$}
\psfrag{l13}[][bl]{$E_q'$} \psfrag{l14}[][bl]{$E_{p+1}'$}
\psfrag{l15}[][bl]{$E_{p+2}'$} \psfrag{l17}[][bl]{$E_q''$}
\psfrag{l18}[][bl]{$I_1\stackrel{\Omega_2^\uparrow}{\longrightarrow}
I_2$} \psfrag{l19}[][bl]{$I_1
\stackrel{\Omega_3}{\longrightarrow}I_2$} \psfrag{v}[bl][bl]{$v$}
\psfrag{P}[bl][bl]{$P$} \psfrag{x}[bl][bl]{$x$}
\psfrag{o1}[B][bl]{$\Omega_1$}
\psfrag{o1u}[B][bl]{$\Omega_1^\uparrow$}
\psfrag{o1d}[B][bl]{$\Omega_1^\downarrow$}
\psfrag{o2}[B][bl]{$\Omega_2$}
\psfrag{o2u}[B][bl]{$\Omega_2^\uparrow$}
\psfrag{o2d}[B][bl]{$\Omega_2^\downarrow$}
\psfrag{o3}[B][bl]{$\Omega_3$} \psfrag{d1}[B][bl]{$D_1$}
\psfrag{d2}[B][bl]{$D_2$} \psfrag{d3}[B][bl]{$D_3$}
\psfrag{d}[B][bl]{$D$} \psfrag{dd}[B][bl]{$D'$}
\psfrag{d1d}[B][bl]{$D_1'$} \psfrag{d2d}[B][bl]{$D_2'$}
\psfrag{d3d}[B][bl]{$D_3'$} \psfrag{dc}[B][bl]{$D \cup C$}
\psfrag{e1}[][bl]{$E_1$} \psfrag{co1}[][bl]{$\Omega_1$}
\psfrag{co1u}[][bl]{$\Omega_1^\uparrow$}
\psfrag{co1d}[][bl]{$\Omega_1^\downarrow$}
\psfrag{co2u}[][bl]{$\Omega_2^\uparrow$}
\psfrag{co2d}[][bl]{$\Omega_2^\downarrow$}
\psfrag{co3}[][bl]{$\Omega_3$}
\psfrag{cso2u}[][bl]{$seq(\Omega_2^\uparrow$)}
\psfrag{ww}[][bl]{$D_1'$} \psfrag{cd1}[][bl]{$D_1$}
\psfrag{cd2}[][bl]{$D_2$} \psfrag{cd3}[][bl]{$D_3$}
\psfrag{cd}[][bl]{$D$} \psfrag{cdd}[][bl]{$D'$}
\psfrag{ca1}[][bl]{$D_2'$} \psfrag{ca2}[][bl]{$D_2'$}
\psfrag{cd3d}[][bl]{$D_3'$} \psfrag{cdc}[][bl]{$D \cup C$}
\psfrag{ce1}[][bl]{$E_1$} \psfrag{dot}[][bl]{\ldots}
\psfrag{e}[][bl]{$=$}
\psfrag{cbd}[][bl]{$E_1'$}\psfrag{w}[][bl]{$E_2'$}
\psfrag{ce2}[][bl]{$E_2$}  \psfrag{e3}[][bl]{$E_3$}
\psfrag{e1d}[][bl]{$E_1'$} \psfrag{e2d}[][bl]{$E_2'$}
\psfrag{e3d}[][bl]{$E_3'$}
\psfrag{a1u}[][bl][2][90]{$\rightsquigarrow$}
\psfrag{a1}[][bl][2][0]{$\rightsquigarrow$} \psfrag{en}[][bl]{$E_n$}
\psfrag{end}[][bl]{$E_n'$} \psfrag{tm}[][bl]{$D_1'$}
\psfrag{tmm}[][bl]{$L_2$}

\begin{document}

\maketitle

\begin{abstract}
We show that, up to ambient
isotopy, the exterior of a hyperbolic knot in $S^3$ admits finitely many bridge punctured 2-spheres of given Euler characteristic, and that there is an algorithm to find all of these surfaces. This yields an algorithm to detect bridge number for hyperbolic knots.
\end{abstract}

\setcounter{page}{1}

\section{Introduction}

\let\thefootnote\relax\footnotetext{MSC (2010): 57M25, 57M27}
Much of knot theory is concerned with understanding knot invariants. Of these
one of the most natural and widely studied is bridge number. However, in common
with other natural knot invariants such as unknotting number, calculating the
bridge number of specific knots can be difficult in practice. The goal of this paper
is to prove, by means of the the following theorem, that for hyperbolic knots it is algorithmically decidable.

\begin{theorem} \label{maintheorem}
Let $K$ be a hyperbolic knot in $S^3$.
Let $M$ be the exterior of $K$ in $S^3$. Then, up to ambient
isotopy, there are only finitely many bridge punctured 2-spheres for
$M$ of given Euler characteristic. Furthermore there is an algorithm
to find all of these surfaces.
\end{theorem}

\begin{corollary} \label{maincorollary}
There exists an algorithm to determine
the bridge number of a hyperbolic knot in $S^3$.
\end{corollary}

The algorithm of Theorem \ref{maintheorem} has as input any diagram for $K$ and an integer $n$, and produces as output a finite list of all
bridge punctured 2-spheres for $M$ of Euler characteristic $n$.
There is no guarantee, however, that the surfaces in the
list are pairwise non-isotopic.

The general scheme that one would
like to follow to prove Theorem \ref{maintheorem} is to search for bridge
punctured 2-spheres by arranging for them to sit in normal or almost
normal form with respect to some triangulation of $M$. These notions were developed, respectively, by Wolfgang Haken in the 1960s, to solve problems such as the recognition problem for Haken 3-manifolds, and more recently by Hyam Rubinstein whose work could be applied in situations where classical normal surface theory breaks down. There are two technical obstructions to applying (almost) normal surface theory to search bridge surfaces of knots, however. Firstly, bridge punctured
2-spheres are, much like Heegaard surfaces, highly compressible in general, and
may well have disjoint compression discs on each side. For this
reason we start by applying work of Marion Moore Campisi \cite{marionalpha}, who builds on work of Chuichiro Hayashi and Koya
Shimokawa \cite{hsh}, to untelescope a bridge punctured 2-sphere
for a knot in $S^3$ to give rise to a generalized bridge surface, something analogous to a generalized Heegaard splitting,
which has more restricted compression discs. The precise result we use is encapsulated in the following proposition, the terms of which are defined in Section 2.

\begin{namedtheorem} For any generalised bridge surface
$\mathcal{B}$ of $M$, a knot exterior, there exists a generalised bridge surface
$\mathcal{B}'$ with the following properties:
\begin{enumerate}
  \item $\mathcal{B}$ may be obtained from $\mathcal{B}'$ by amalgamation.
  \item Each thin surface of $\mathcal{B}'$ is incompressible.
  \item Each thick surface  $K_i$ of $\mathcal{B}'$ is strongly
irreducible in $M_i$, the submanifold of $M$ obtained by cutting along the adjacent thin surfaces $N_i$ and $N_{i-1}$.
  \item No thick surface $K_i$ cobounds a product $(\textrm{Surface}) \times
  I$ with an adjacent thin surface $N_i$ or $N_{i-1}$.
\end{enumerate}
\end{namedtheorem}

A second, more fundamental, obstruction to applying  normal surface theory to bridge surfaces relates to the algorithmic side of the theory. As exposed by Sergei Matveev in \cite{mat}, it is essential to gain control, one way or another, of normal surfaces of non-negative Euler characteristic. Moreover, if one works in an arbitrary triangulation of a knot exterior, normal tori are endemic. While the prevalence of normal 2-spheres can be effectively controlled using 0-efficient triangulations, developed by William Jaco and Hyam Rubinstein \cite{normtori}, efforts to control normal tori via a corresponding theory of 1-efficient triangulations have so far proved more difficult. An approach to this problem for the case where thin position coincides with bridge position using 1-efficient triangulations is given in work of David Bachman \cite{bach}. Robin Wilson has also proved some related non-algorithmic results about bridge surfaces in triangulated knot exteriors in \cite{robinbridge}, as have David Bachman, Ryan Derby-Talbot and Eric Sedgwick in a  more general setting \cite{bachderbsedg}.

The approach we use in this paper is different: We work in decompositions of hyperbolic knot exteriors called partially flat angled ideal triangulations. Their existence and use is summarized, respectively, by the following two theorems of Marc Lackenby \cite{tn1alg}:

\begin{namedtheoremc} [Theorem 2.2 of \cite{tn1alg}] Any finite-volume hyperbolic $3$-manifold $M$ with
non-empty boundary has a partially flat angled ideal triangulation.
Moreover, there is an algorithm that constructs one, starting with
any triangulation of
$M$.
\end{namedtheoremc}

\begin{namedtheoremd} [Theorem 2.1 of \cite{tn1alg}] Let $T$ be a partially flat angled ideal
triangulation of a $3$-manifold $M$. Then any connected $2$-normal surface in $T$
with non-negative Euler characteristic is normally parallel to a
boundary
component.
\end{namedtheoremd}

Partially flat angled ideal triangulations are defined in \cite{tn1alg}, and in Section 4 of this paper. They are ideal triangulations decorated with some extra combinatorial data, often arising naturally from to the  geometry of a totally geodesic
ideal polyhedral decomposition. See \cite{ep} and \cite{tn1alg}. The 2-normal surfaces in the later of the two theorems above are surfaces which consist of normal triangles, quadrilaterals or octagons. Thus normal tori are controlled. Now, in  \cite{tn1alg} Lackenby was concerned with almost normal Heegaard surfaces in partially flat angled ideal 
triangulations. Since Heegaard surfaces are closed, it suffices in that context to consider genuine ideal triangulations, that is decompositions for non-compact $3$-manifolds built out of tetrahedra with their vertices removed. In this paper we consider surfaces that run up to the toral boundary of a knot exterior, so the decompositions we consider are slightly different, being built out of truncated tetrahedra, that is tetrahedra with a small open neighborhood of their vertices removed, rather than ideal tetrahedra. Nevertheless, in order to preserve some unity of language between  \cite{tn1alg} and this paper, we will refer to these decompositions for compact manifolds as ideal triangulations also. Now, a truncated tetrahedron is a polyhedron with eight faces and
eighteen edges. This means that the possible configurations of compact
surfaces in  ideally triangulated compact 3-manifolds are somewhat
complicated. In Section 3, following arguments which are inspired (more philosophically than in application) by Michelle Stocking's treatment of almost normal Heegaard splittings \cite{stocking}, we prove that the generalized bridge surfaces of interest, namely those arising as $\mathcal{B}'$ in Theorem 2.1, can be arranged
to sit `nicely' within our ideally triangulated knot exterior $M$.

This brings us on to the question of what we mean by `nicely'. It is perhaps surprising that the change in polyhedral decomposition, from triangulation to ideal triangulation, means that many of the arguments that can be applied in the context of triangulations break down when one attempts to apply them to ideally triangulated compact 3-manifolds. For example, with most natural definitions of what `normal' should mean in such an ideally triangulated 3-manifold, it fails to be the case that normal surfaces are incompressible in the complement of the 1-skeleton. Another complication is that when we untelescope a bridge surface via Theorem 2.1, our notion of strong irreducibility that is held by the resulting surface $\mathcal{B}'$ fails to provide any control on boundary compression discs. To escape from this quagmire we adopt a goal-oriented approach, relaxing the requirements of what `nicely' should mean somewhat and introducing the notions of interiorly-normal, interiorly-normal to one
side and almost interiorly-normal surfaces which are considerably more general than normal and almost normal surfaces. The following theorem, whose assumptions are chosen to tie in closely with Theorem 2.1, is proved in Section 3.

\begin{namedtheoreme} Let $K$ be a knot is $S^3$. Let $M$ be the exterior of $K$ and let $M$ have an ideal triangulation which contains no 2-spheres consisting of just triangles, squares and octagons. Let $\mathcal{B}$ be a generalised bridge surface for $M$ such that:
\begin{enumerate}
  \item Each thin surface of $\mathcal{B}$ is incompressible.
  \item Each thick surface $K_i$ of $\mathcal{B}$  is strongly
irreducible in $M_i$, the submanifold of $M$ obtained by cutting along the adjacent thin surfaces $N_i$ and $N_{i-1}$.
  \item No thick surface $K_i$ cobounds a product $(\textrm{Surface}) \times
  I$ with an adjacent thin surface $N_i$ or $N_{i-1}$.
\end{enumerate}
Then $\mathcal{B}$ may be ambient isotoped so that:
\begin{enumerate}
  \item The thin surfaces of $\mathcal{B}$ are interiorly normal.
  \item The thick surfaces of $\mathcal{B}$ are interiorly-normal, interiorly-normal to one
side or almost interiorly-normal.
  \end{enumerate}
\end{namedtheoreme}

The definitions of interiorly-normal, interiorly-normal to one
side and almost interiorly-normal surfaces are given in Section 2. They are somewhat broader classes of surface than is usually considered when generalizing (almost) normal surface theory. For example, there are infinitely many interiorly-normal `disc types' in a truncated tetrahedron, in contrast to the seven disc types per tetrahedron that arise in classical normal surface theory. This flexibility needs to be compensated for by gaining more precise control of how a generalized bridge surface intersects the boundary of $M$. It is with this in mind that we prove in Section 4 the following suitably empowered variation of Theorem 3.1. 

\begin{namedtheoremg} Let $K$ be a knot is $S^3$. Let $M$ be the exterior of $K$ and let $M$ have an ideal triangulation which contains no 2-spheres consisting of just triangles, squares and octagons. Let $F$ be bridge punctured $2$-sphere for $M$. Let $\mathcal{B'}$ be a generalized bridge surface for $M$ such that:
\begin{enumerate}
  \item Each thin surface of $\mathcal{B}'$ is incompressible.
  \item Each thick surface $K_i$ of $\mathcal{B'}$ is strongly
irreducible in $M_i$, the submanifold of $M$ obtained by cutting along the adjacent thin surfaces $N_i$ and $N_{i-1}$.
  \item No thick surface $K_i$ cobounds a product $(\textrm{Surface}) \times
  I$ with an adjacent thin surface $N_i$ or $N_{i-1}$.
  \item  $\mathcal{B}'$ amalgamates to give $F$.
\end{enumerate}
Then there is a computable constant $b$, which may be calculated
from the Euler characteristic of $F$ and the ideal triangulation
of $M$, such that $\mathcal{B'}$ may be ambient isotoped so that:
\begin{enumerate}
  \item The thin surfaces of $\mathcal{B'}$ are interiorly normal.
  \item The thick surfaces of $\mathcal{B'}$ are interiorly-normal, interiorly-normal to one
side or almost interiorly-normal.
\item $\mathcal{B'}$ intersects the boundary 1-skeleton of $T$ at most $b$ times.
  \end{enumerate}
\end{namedtheoremg}

The main technical algorithmic theorem of this paper is the following. It is this theorem that brings into focus the usefulness of interiorly-normal, interiorly-normal to one side and almost interiorly-normal surfaces. Its proof appears in Section 4 of this paper.

\begin{namedtheoremf} Let $T$ be a partially flat angled ideal
triangulation of a 3-manifold $M$. Then, for any positive integers $n$ and $b$,
$T$ contains only finitely many properly embedded surfaces $F$ such
that:
\begin{enumerate}
  \item $- \chi(F) \leq n$.
  \item $F$ intersects the boundary 1-skeleton of $T$ at most $b$ times.
  \item Each component of $F$ is either interiorly-normal,
  interiorly-normal to one side or almost interiorly-normal.
\end{enumerate}
Furthermore, there is an algorithm to construct each of these
surfaces.
\end{namedtheoremf}

This paper is arranged as follows: Section 2 is preliminaries. In Section 3 we show
that each component of a generalised bridge surface of the type we
wish to search for may be ambient isotoped so that it is
interiorly-normal, interiorly-normal to one side or almost
interiorly-normal. In Section 4 we focus on the algorithmic tools needed for Theorem \ref{maintheorem} and in particular we prove Theorems 4.3 and 4.6. We also comment in Section 4 on how the process of amalgamation is achieved algorithmically and bring the various results together to show how to construct
the algorithm of Theorem \ref{maintheorem}. 

The author is grateful to Marc Lackenby for many helpful discussions during the preparation of this paper. 

\section{Preliminaries}

\begin{figure}[h!]
\centering
\includegraphics{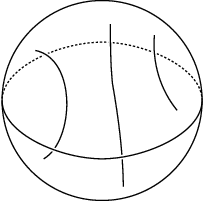}
\caption{A trivial tangle} \label{trivtang}
\end{figure}

A \emph{tangle} is a pair $(B,T)$ where $B$ is a 3-ball and $T$ is a
finite collection of disjoint arcs properly embedded in $B$. A
\emph{trivial tangle} is a tangle that is homeomorphic as a pair to
$(D \times I,P \times I)$ where $D$ is a disc, $I$ is the closed
unit interval, and $P$ is a finite collection of points in the
interior of $D$. Figure \ref{trivtang} shows an example of a trivial
tangle.

Let $K \subseteq S^3$ be a knot. Let $F$ be a 2-sphere in $S^3$
that satisfies the following properties:
\begin{enumerate}
\item $K$ intersects $F$ transversely.
\item $F$ cuts $(S^3,K)$ into two components, both of which are
trivial tangles.
\end{enumerate}
Then $F$ is known as a \emph{bridge 2-sphere} for $K$. The minimum
of $\frac{|F \cap K |}{2}$ over all bridge 2-spheres, $F$, is known
as the \emph{bridge number} of $K$.

If we remove a small open neighbourhood of $K$ from $S^3$, then a
bridge 2-sphere becomes a \emph{bridge punctured 2-sphere} for the
knot exterior, and the trivial tangles become \emph{trivially
punctured 3-balls}. The part of a trivially punctured 3-ball which
coincides with the bridge punctured 2-sphere is called the
\emph{outside boundary} and the closure of the rest of the boundary is called the
\emph{inside boundary}. In general, if $N$ is a 3-dimensional
submanifold of the exterior, $M$, of a knot in $S^3$, then we shall
call the closure of that part of $\partial N$ which is disjoint from
$\partial M$ the \emph{outside boundary} of $N$, and we shall denote
this by $\partial^* N$.

\begin{figure}[h!]
\centering
\includegraphics{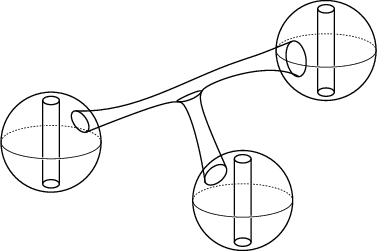}
\caption{A trivially punctured 3-ball} \label{trivpunct3ball}
\end{figure}

Let $M$ be the  3-manifold with boundary obtained by
removing a small open neighbourhood of a  knot from $S^3$.
Let $F$ be a bridge punctured 2-sphere for $M$. Then $M$ decomposes,
when cut along $F$, into two trivially punctured 3-balls, $B_1$ and
$B_2$. Consider a small regular neighbourhood of a meridian curve of
one of the punctures. This is a trivially punctured 3-ball with only
one puncture. We shall call such an object a \emph{punctured
0-handle}. A trivially punctured 3-ball, $B$, may be constructed by
taking one punctured 0-handle around a meridian of each puncture of
$B$ and connecting them with 1-handles, as shown in Figure
\ref{trivpunct3ball}. In other words, there is a way of constructing
a trivially punctured 3-ball by starting with a collection of
punctured 0-handles and attaching a collection of 1-handles.

This construction may be applied to both $B_1$ and $B_2$. However,
in a similar way to the way one uses a Heegaard splitting to
determine a handle decomposition for a closed 3-manifold, we will
refer to the 1-handles in $B_2$ as 2-handles and the punctured
0-handles as \emph{punctured 3-handles}. To recap, we have built $M$
by starting with a collection of punctured 0-handles, attaching a
collection of 1-handles, then a collection of 2-handles, and finally
a collection of punctured 3-handles. We shall call such a
construction a \emph{bridge decomposition} of $M$.

With the same ideas in mind as in \cite{untel}, suppose that we were
to build $M$ by starting with a collection of punctured 0-handles,
$H^0$, then attaching a collection of 1-handles, $H^1_1$, then a
collection of 2-handles, $H^2_1$, and then some more 1-handles,
$H^1_2$, and some more 2-handles, $H^2_2$, etc ... and then some
more 1-handles, $H^1_n$, then a collection of 2-handles, $H^2_n$,
and finally a collection of punctured 3-handles, $H^3$. We shall
refer to such a construction as a \emph{generalised bridge
decomposition} of $M$.

Let $N_0 = \partial^* H^0$, and for $i=1, \ldots ,n$ let
$$N_i =  \partial^* (H^0 \cup \bigcup_{k=1}^i ({H^1_k \cup H^2_k})).$$

Let $K_1 = \partial^* (H^0 \cup H^1_1)$ and for $i=2, \ldots ,n$ let
$$K_i = \partial^* (H^0 \cup H^1_1 \cup \bigcup_{k=2}^i ({H^1_k \cup
H^2_{k-1}})).$$

Note that the surfaces $N_i$ and $K_i$ defined in this way are not
disjoint. Carry out a small isotopy to rectify this. We shall refer
to the resulting surfaces $N_i$ (resp. $K_i$) as the thin (resp. thick)
surfaces of the decomposition. Let $\mathcal{B}$ denote the
collection of surfaces $N_i$ and $K_i$. We shall refer to
$\mathcal{B}$ as a \emph{generalised bridge surface} for $M$.

Let $M_i$ be the submanifold of $M$ whose boundary consists of
$N_{i-1}$ and $N_i$, so that $K_i$ lies inside $M_i$. We will say
that $K_i$ is \emph{strongly irreducible} in $M_i$ if any two
compression discs for $K_i$ in $M_i$ on opposite sides of $K_i$
intersect at some point along their boundary.

In the event that one of the surfaces $K_i$ is not strongly
irreducible, we may perform an \emph{untelescoping operation} on the
generalised bridge decomposition. Untelescoping operations only
affect the 1-handles and 2-handles of the decomposition, and they
are described in \cite{untel}. The reverse procedure of an
untelescoping operation is called \emph{amalgamation}. For an informative
description of amalgamation see \cite{heeggenuslack}. Note that the
use of punctured 0-handles and punctured 3-handles instead of
0-handles and 3-handles does not affect these notions. Further note
that amalgamation may be carried out algorithmically, as described
in \cite{tn1alg}.

The following thoerem is absolutely key in the proof of Theorem 1.1:

\begin{proposition} \label{untelprop} For any generalised bridge surface
$\mathcal{B}$ of $M$, a knot exterior, there exists a generalised bridge surface
$\mathcal{B}'$ with the following properties:
\begin{enumerate}
  \item $\mathcal{B}$ may be obtained from $\mathcal{B}'$ by amalgamation.
  \item Each thin surface of $\mathcal{B}'$ is incompressible.
    \item Each thick surface $K_i$ of $\mathcal{B}'$ is strongly
irreducible in $M_i$, the submanifold of $M$ obtained by cutting along the adjacent thin surfaces $N_i$ and $N_{i-1}$.
  \item No thick surface $K_i$ cobounds a product $(\textrm{Surface}) \times
  I$ with an adjacent thin surface $N_i$ or $N_{i-1}$.

\end{enumerate}
\end{proposition}

For a proof of Theorem 2.1 the reader is referred to Marion Moore Campisi's paper \cite{marionalpha}. Note that in that paper  Campisi refers to 0-beads and 2-beads instead of punctured 0-handles and punctured 3-handles and that her result applies in somewhat more generality.

Our strategy for
algorithmically searching for bridge surfaces, $\mathcal{B}$, will
be to search for generalised bridge decompositions, $\mathcal{B}'$,
as in Theorem 2.1, and then algorithmically amalgamate
$\mathcal{B}'$. This is achieved by finding an ideal triangulation
together with some extra structure for $M$ and placing the thin
surfaces of $\mathcal{B}'$ into something which resembles normal
form and the thick surfaces of $\mathcal{B}'$ into something which
resembles almost normal form.

\begin{figure}[h!]
\centering
\includegraphics{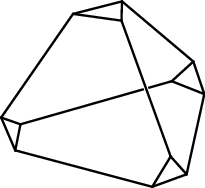}
\caption{A truncated tetrahedron} \label{cit}
\end{figure}

A \emph{truncated tetrahedron} is a tetrahedron with a
small open neighbourhood of its vertices removed, as shown in Figure
\ref{cit}. We will occasionally abuse language slightly and  refer to a truncated tetrahedron as an ideal tetrahedron or simply a tetrahedron.

A truncated tetrahedron  has eight faces. Four of these are
triangular and the other four are hexagonal. These faces are called
\emph{boundary faces} and \emph{interior faces} respectively. If one
forms a 3-manifold by  homeomorphically identifying in pairs the
interior faces of a collection of truncated tetrahedra so that the boundary faces patch together to form a surface, then
the 3-manifold will be said to have a \emph{compact ideal
triangulation}. We will sometimes refer to a compact ideal triangulation as an ideal triangulation or simply a triangulation when there is no possibility of confusion. The edges of the truncated tetrahedra manifested in the
resulting 3-manifold are of two types, namely those that which lie
on the boundary of a boundary face and those that do not. These are
called \emph{boundary edges} and \emph{interior edges} respectively.

A properly embedded arc on a boundary face
of a truncated tetrahedron is said to be a \emph{normal arc} if it has
endpoints on different boundary edges. A properly embedded arc on an
interior face is said to be an \emph{interiorly-normal arc} if it
has endpoints on different edges or the same boundary edge.

\begin{figure}[h!]
\centering
\includegraphics{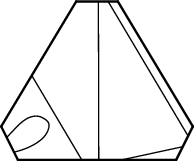}
\caption{Examples of interiorly-normal arcs} \label{intnormarcs}
\end{figure}

Examples of different types of interiorly-normal arc are shown in
Figure \ref{intnormarcs}.

\begin{figure}[h!]
\centering
\includegraphics{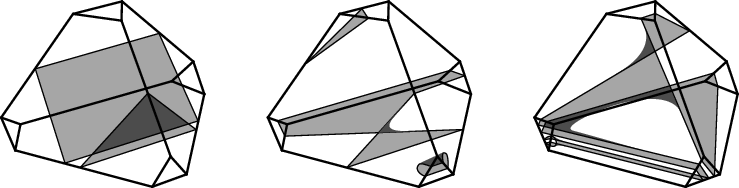}
\caption{Some interiorly normal discs} \label{exintdisc}
\end{figure}

\begin{define}Let $T$ be a truncated tetrahedron.
Let $D \subseteq T$ be a properly embedded disc in general position
with respect to the 1-skeleton which satisfies the following
conditions:
\begin{enumerate}
\item $\partial D$ intersects each boundary face of $T$ in normal
arcs.
\item $\partial D$ intersects each interior face of $T$ in interiorly-normal arcs.
\item $\partial D$ intersects each interior edge at most once.
\end{enumerate}
Then $D$ is said to be an \emph{interiorly-normal disc}. Let $M$ be
a 3-manifold with boundary with compact ideal triangulation. A
properly embedded surface $S \subseteq M$ is said to be an
\emph{interiorly-normal surface} if it intersects each truncated  tetrahedron in interiorly-normal discs.
\end{define}
Some examples of interiorly-normal disc are shown Figure
\ref{exintdisc}.

\begin{define}Let $S \subseteq M$ be an embedded surface
in general position with respect to the 1-skeleton of an ideally triangulated compact 3-manifold. Let $E$ be a
disc embedded in $T$, a truncated tetrahedron in an ideal
triangulation of $M$, whose interior is disjoint from $S \cup
\partial T$  and $\partial E = \alpha \cup \beta$ where $\alpha \cap
\beta = \partial \alpha =
\partial \beta$, $\alpha$ is an arc in $S \cap T$ and $\beta$ is a
sub-arc of an interior-edge of $T$. Then $E$ is said to be an
\emph{edge compression disc} for $S$. If $\beta$ is instead an arc embedded
in the interior of an interior face of $T$ then we shall say that
$E$ is a \emph{face compression disc} for $S$.
\end{define}

\begin{define}Let $T$ be a truncated tetrahedron.
Let $D \subseteq T$ be a properly embedded disc in general position
with respect to the 1-skeleton which satisfies the following
conditions:
\begin{enumerate}
\item $\partial D$ intersects each boundary face of $T$ in normal
arcs.
\item $\partial D$ intersects each interior face of $T$ in interiorly-normal arcs.
\item $D$ admits at least one edge compression disc.
\item All edge compression discs emanate from the same side of $D$.
\end{enumerate}
Then $D$ is said to be \emph{interiorly-normal to one side}. A
2-sided properly embedded surface $S \subseteq M$ in general
position with respect to the 1-skeleton is said to be
\emph{interiorly-normal to one side} if it intersects each
truncated tetrahedron in interiorly-normal and
interiorly-normal to one side discs, at least one component of
intersection of $S$ with a truncated tetrahedron admits an edge
compression disc and all such edge compression discs emanate from
the same side of $S$. We shall refer to the side without edge
compression discs as the \emph{interiorly-normal side}.
\end{define}

\begin{remark} We will sometimes say that a surface is
interiorly-normal on a particular side. When we do, we include the
possibility of the surface being interiorly-normal.
\end{remark}

\begin{define}Let $T$ be a truncated tetrahedron.
Let $D \subseteq T$ be a properly embedded disc in general position
with respect to the 1-skeleton which satisfies the following
conditions:
\begin{enumerate}
\item $\partial D$ intersects each boundary face of $T$ in normal
arcs.
\item $\partial D$ intersects each interior face of $T$ in interiorly-normal arcs.
\item $\partial D$ intersects each interior edge at most twice.
\item $D$ admits at least one edge compression disc on each side.
\item Any pair of edge compression discs
for $D$ emanating from opposite sides of $D$ intersect.
\end{enumerate}
Then $D$ is said to be an \emph{almost interiorly-normal disc}.
\end{define}
\begin{define}Let $T$ be a truncated tetrahedron.
Let $A \subseteq T$ be a properly embedded annulus in general
position with respect to the 1-skeleton which is formed by
connecting two interiorly-normal discs with a tube parallel to a
face of the interior 2-skeleton. Then $A$ is said to be an
\emph{almost interiorly-normal annulus}.
\end{define}
\begin{define}A properly embedded surface $S \subseteq M$
in general position with respect to the 1-skeleton is said to be
\emph{almost interiorly-normal} if it intersects each truncated
tetrahedron in interiorly-normal discs, apart from in
precisely one tetrahedron which it intersects in interiorly-normal
discs and precisely one almost interiorly-normal disc or annulus.
\end{define}
\begin{remark} An interiorly-normal (resp. almost
interiorly-normal) surface which does not intersect the boundary of
$M$ is normal (resp. almost normal) in the
classical sense. See \cite{stocking}.
\end{remark}
Note that there are infinitely many interiorly-normal disc types in
a given truncated tetrahedron. Later on we will need to restrict this
class of discs to a finite collection. The following definition
shall be key in this respect.

\begin{define} The boundary (resp. interior) edge degree of
a properly embedded surface $S \subseteq M$ is the number of
intersections of $S$ with the boundary (resp. interior) 1-skeleton.
\end{define}
\begin{define} An arc embedded in a face of the interior
2-skeleton will be said to be a \emph{normal arc} if it joins two
different interior edges of the interior face on which it lies. This
agrees with how normal arcs are usually defined in non-ideal
tetrahedra. A \emph{normal curve} is a curve embedded on the
boundary of a truncated tetrahedron which consists of normal arcs. The
\emph{length} of a normal curve is the number of normal arcs which
it consists of.
\end{define}

\begin{define} Let $S \subseteq M$ be a properly embedded
surface. Let $C \subseteq \textrm{int}(M)$ be an embedded arc such
that $C \cap S = \partial C$. In the event that $M$ has an ideal
triangulation suppose that $C$ does not intersect the interior
1-skeleton and that $C$ intersects the interior 2-skeleton
transversely in a finite number of points. Let $D$ be a disc and let
$C \times D$ be a small product neighbourhood of $C$ such that $(C
\times D) \cap S = (C \cap S) \times D$. Define a new surface $$S' =
(S \cup (C \times
\partial D))\backslash (\partial C \times \textrm{int}(D)).$$ We
shall say that \emph{$S'$ is obtained from $S$ by adding a tube
along $C$}.

Let $G$ be an embedded graph in $\textrm{int}(M)$ with at least one
1-valent vertex in each connected component. Suppose that each
1-valent vertex lies on $S$ and that $S$ does not intersect $G$
other than at 1-valent vertices. In the event that $M$ has an ideal
triangulation suppose that $G$ does not intersect the interior
1-skeleton and that $G$ intersects the interior 2-skeleton
transversely in a finite number of points, none of which are
vertices of $G$. Place a small 2-sphere at each vertex of $G$ with
valence more than 1. Now attach a tube along each edge of $G$ and
call the resulting surface $S'$. In this case we shall say that
\emph{$S'$ is obtained from $S$ by adding tubes along $G$}. We shall
refer to $G$ as the \emph{core of the tubes of $S'$}.

If $S$ and $S'$ are isotopic then we shall say that the tubes are
\emph{non-essential}. Otherwise they are \emph{essential}.
\end{define} 
\section{Isotoping generalised bridge surfaces}

The entirety of this section is devoted to proving the following theorem.

\begin{theorem} \label{nktheorem} Let $K$ be a knot is $S^3$. Let $M$ be the exterior of $K$ and let $M$ have an ideal triangulation which contains no 2-spheres consisting of just triangles, squares and octagons. Let $\mathcal{B}$ be a generalised bridge surface for $M$ such that:
\begin{enumerate}
  \item Each thin surface of $\mathcal{B}$ is incompressible.
  \item Each thick surface $K_i$ of $\mathcal{B}$ is strongly
irreducible in $M_i$, the submanifold of $M$ obtained by cutting along the adjacent thin surfaces $N_i$ and $N_{i-1}$.
  \item No thick surface $K_i$ cobounds a product $(\textrm{Surface}) \times
  I$ with an adjacent thin surface $N_i$ or $N_{i-1}$.
\end{enumerate}
Then $\mathcal{B}$ may be ambient isotoped so that:
\begin{enumerate}
  \item The thin surfaces of $\mathcal{B}$ are interiorly normal.
  \item The thick surfaces of $\mathcal{B}$ are interiorly-normal, interiorly-normal to one
side or almost interiorly-normal.
  \end{enumerate}
\end{theorem}

Our first objective is to arrange the boundaries of the surfaces $N_i$ and $K_i$ on $\partial M$. With this in mind, note that the boundary torus of $M$ admits a natural product structure,
$\partial M = S^1_M \times S^1_L$, where the first factor denotes
the meridional coordinate, the second factor denotes the
longitudinal coordinate and the boundary circles of the surfaces
$N_i$ and $K_i$ are constant on the longitudinal factor. Let $T$ be
the compact ideal triangulation of $M$. We shall start by
isotoping $T$ so that the boundary 1-skeleton satisfies the
following properties:
\begin{enumerate}
  \item All the boundary edges are transverse to the foliation of $\partial
  M$ by meridional circles.
  \item All the vertices of the triangulation have different
  meridional coordinates.
\end{enumerate}
We may find an isotopy to satisfy the first property by
\cite{straighttorus}, in which it is proved that any simple
triangulation (see \cite{straighttorus}) of a torus may be isotoped
so that all the edges are geodesics in the Euclidean metric. Note
that a simple Euler characteristic calculation implies that the
boundary 1-skeleton is a simple triangulation of $\partial M$. If
necessary, a small isotopy suffices to ensure that the second
property holds.

Next proceed by isotoping the vertices of the triangulation into
$H_0$ as follows. Let $f$ be a homeomorphism of $S^1$ isotopic to
the identity. Let $f_t:S^1 \rightarrow S^1$ be an isotopy of $S^1$
so that $f_1 = f$ and $f_0$ is the identity on $S^1$. Let $\partial
M \times [0,1]$ be a small collar of $\partial M$ where $\partial M
\times \{0\} =
\partial M$. Define $F:M \rightarrow M$ by $F(x) = x$ when $x
\notin \partial M \times [0,1]$ and for $(m,l,s) \in S^1_M \times
S^1_L \times [0,1]$ by $$F(m,l,s) = (m,f_{(1-s)}(l),s).$$ We shall
call an isotopy constructed in this way a \emph{boundary
height adjusting isotopy}.

\begin{figure}[h!]
\centering
\includegraphics{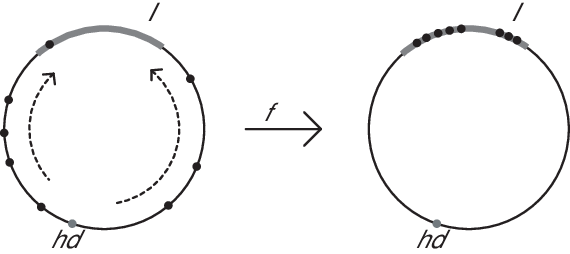}
\caption{An isotopy of $S_L^2$} \label{circlehomeo}
\end{figure}

Let $h_1 \ldots h_l \in S^1_L$ denote the longitudinal coordinates
of the vertices of the ideal triangulation. Let $I \in S^1_L$ be a
subinterval of $S^1_L$ for which $(S^1_M \times I) \subseteq H_0$.
Let $h' \in S_L^1 \backslash (I \cup \bigcup_{i=1}^l{\{h_i\}})$. Now
define $f:S^1_L  \rightarrow S^1_L$ so that $f$ fixes $h'$ and
$f(h_i) \in I$ for $i = 1 \ldots l$, as shown in Figure
\ref{circlehomeo}.

Note that $f$ is isotopic to the identity. Let $f_t$ be an isotopy
from the identity to $f$. Then $F$, the boundary height adjusting
isotopy associated to $f_t$, is an isotopy of $M$ which sends all
the vertices of $T$ into $H_0$ and furthermore keeps the boundary
edges transverse to the foliation of $\partial
  M$ by meridional circles.

It is worth emphasising that the following properties still hold:
\begin{enumerate}
  \item The boundary edges of the triangulation are transverse to the foliation of $\partial
  M$ by meridional circles.
  \item The boundary circles of the surfaces
$N_i$ and $K_i$ have constant longitudinal coordinate on $\partial
M$.
\end{enumerate}
Together these two properties mean that the boundary circles of the
surfaces $N_i$ and $K_i$ intersect the boundary faces of the
triangulation in normal arcs. Later on we will show how to obtain a
bound on the boundary edge degrees of these surfaces, but for now the
next step is to isotope the surfaces $N_i$ rel boundary into
interiorly-normal position with respect to the triangulation. This
we achieve with the following proposition.

\begin{proposition} \label{intnormprop}Let $M$ be a compact ideally triangulated
3-manifold with boundary. Let $S \subseteq M$ be a properly embedded
surface whose boundary intersects the boundary faces of the ideal
triangulation in normal arcs. Then $S$ may be ambient isotoped rel
boundary, compressed and have 2-sphere components removed to lie in
interiorly-normal position with respect to the ideal triangulation.
\end{proposition}

\noindent \textbf{Proof.} We shall use similar techniques to those used to
prove Theorem 3.3.21 in \cite{mat}. Start by isotoping $S$ rel
boundary into general position with respect to the triangulation. We
will need some different measures of complexity for $S$. Recall the
\emph{interior edge degree} of $S$, $e(S)$, is the total number of
intersections of $S$ with the interior 1-skeleton of $M$. Let
$\gamma(S) = \sum_{i=1}^{m}(1-\chi(S_i))$ where $S_1,\ldots,S_m$ are
the connected components of the intersection of $S$ with each
tetrahedron which are not 2-spheres. Let $n(S)$ be the total number
of connected components of $S$. Our measure of complexity for $S$
will be the \emph{weight} of S, $w(S) := (e(S),\gamma(S),n(S)) \in
(\mathbb{N} \cup 0) \times (\mathbb{N} \cup 0) \times (\mathbb{N}
\cup 0)$ where $(\mathbb{N} \cup 0) \times (\mathbb{N} \cup 0)
\times (\mathbb{N} \cup 0)$ is ordered lexicographically. Our
strategy will be to carry out a sequence of moves which all reduce
the weight of $S$. These moves are as follows.
\begin{enumerate}
  \item[$\mathcal{N}_1$] Suppose that the intersection of $S$ with a truncated tetrahedron, $T$, admits a compression disc $D \subseteq
  T$. Then compress $S$ along $D$.
  \item[$\mathcal{N}_2$] Suppose that $S$ admits an edge compression disc,
  $D$. Isotope $S$ across $D$.
  \item[$\mathcal{N}_3$] Suppose that a component of intersection of $S$ with
  a tetrahedron of the triangulation is a 2-sphere. Remove this
  component.
  \item[$\mathcal{N}_4$] Suppose that a component of $S$ is a 2-sphere that
  intersects the interior 2-skeleton of the triangulation in a single circle
  contained in an interior face of the triangulation. Remove this
  component.
\end{enumerate}

Clearly, $\mathcal{N}_2$ decreases $e(S)$ while $\mathcal{N}_1$, $\mathcal{N}_3$ and $\mathcal{N}_4$
preserve $e(S)$. Also $\mathcal{N}_1$ decreases $\gamma(S)$ (see \cite{mat}), while
$\mathcal{N}_3$ and $\mathcal{N}_4$ preserve $\gamma(S)$. Finally, $\mathcal{N}_3$ and $\mathcal{N}_4$ both
reduce $n(S)$. Hence all four moves decrease the weight of $S$.
Hence, after applying these moves as much as possible to $S$ we must
be left with a new surface, $S'$, which does not admit any of the
moves $\mathcal{N}_1$ to $\mathcal{N}_4$, and is obtained from $S$ by means of
compressing, isotoping rel boundary and removing 2-spheres. It
remains to show that $S'$ is interiorly-normal. First note that $S'$
is in general position relative to the ideal triangulation. Also,
$S'$ intersects each truncated tetrahedron in a collection of discs. The
boundary of $S'$ intersects the boundary faces of the ideal
tetrahedra in normal arcs since $S$ does, and so the first
requirement of interior-normality holds. If $S'$ intersects an
interior face in a simple closed curve then we may apply move $\mathcal{N}_1$
or $\mathcal{N}_4$. Hence $S'$ intersects the interior faces in embedded arcs
and the second condition holds. Finally, if one of the discs of
intersection of $S'$ with a truncated tetrahedron has boundary which
intersects an interior edge more than once, then an innermost
curve/outermost arc argument implies that we may apply an $\mathcal{N}_2$
move. Hence the third requirement is satisfied and $S'$ is
interiorly-normal. \hfill $\square$\\

Proposition \ref{intnormprop} immediately yields the following corollary: 

\begin{cor} Let $M$ be a compact irreducible ideally
triangulated 3-manifold with boundary. Let $S \subseteq M$ be a
properly embedded, incompressible surface with no 2-sphere
components which intersects the boundary faces of the ideal
triangulation in normal arcs. Then $S$ may be ambient isotoped rel
boundary into interiorly-normal position with
respect to the ideal triangulation.
\end{cor}

Now, all the surfaces $N_i$  satisfy the hypotheses of Corollary 3.3
and hence they may be isotoped rel boundary into interiorly-normal
position.

Recall that, for $i = 1 \ldots n$, $M_i$ is the sub-manifold of $M$
whose boundary consists of the interiorly-normal surfaces $N_{i-1}$
and $N_i$ as well as part of the boundary of $M$, so that $K_i$ lies
inside $M_i$ for each $i$. Now, $K_i$ is obtained from $N_{i-1}$
(resp. $N_{i}$) by adding tubes. Let $C_1$ (resp. $C_2$) denote the
core of these tubes. Define $h:M_i \rightarrow [0,1]$ by $h(N_{i-1}
\cup C_1) = 0$, $h(N_{i} \cup C_2) = 1$ and note that $M_i
\backslash ((N_{i-1} \cup C_1) \cup (N_{i} \cup C_2)) = K_i \times
(0,1)$ so that for $p \in M_i \backslash ((N_{i-1} \cup C_1) \cup
(N_{i} \cup C_2))$ we may define $h(p)$ by projection onto the
second factor. For a point $p \in M_i$ we shall refer to $h(p)$ as
the \emph{height} of $p$. A surface of the form $h^{-1}(t)$ for $t
\in (0,1)$ shall be called an \emph{interior fibre} of $h$.

Note that, for $i = 1 \ldots n$, $M_i$ contains no vertices of the
ideal triangulation because, before isotoping the surfaces $N_i$
into interiorly-normal form, we apply a suitable boundary height
adjusting isotopy to push all the vertices into $H_0$.

The notion of \emph{thin position} was first introduced by David
Gabai in \cite{gabaithin} with reference to knots, but since then
the notion has found applications in a variety of areas. We will use
an adapted version of thin position here. Consider the interior
1-skeleton of the ideal triangulation in $M_i$ and apply a small isotopy so that $h$ restricts to a Morse function on this. Then as we follow the
fibres of $h$ from $h^{-1}(1)$ down to $h^{-1}(0)$ we see a sequence
of maxima and minima of the interior 1-skeleton. A small isotopy
ensures that there are finitely many and that they occur at
different heights. The levels which intersect one of these maxima or
minima shall be referred to as \emph{critical levels}. The levels in
between two critical levels all look the same relative to the
1-skeleton, and these shall be referred to as \emph{non-critical
levels}. Let $f_1=h^{-1}(a_1),\ldots,h^{-1}(a_m)=f_m$ be
representatives of each non-critical level. Note that each $f_i$
intersects the interior 1-skeleton transversely and define the
\emph{width}, $w(f_i)$, of each non-critical level to be the number
of intersections of $f_i$ with the interior 1-skeleton. Note that if
we isotope the interior 1-skeleton about rel boundary then we may
affect the widths of the non-critical levels. The interior
1-skeleton is said to be in \emph{thin position} with respect to $h$
if the sum of the widths of the non-critical levels is minimal up to
ambient isotopy rel $\partial M_i$ of the interior 1-skeleton. The
sum of the widths of the non-critical levels of the 1-skeleton when
it is in thin position is known as the \emph{width of the interior
1-skeleton} with respect to $h$. A non-critical level that lies
immediately above a maximum but immediately below a minimum is said
to be a \emph{thin level}. One that lies immediately above a
minimum but immediately below a maximum is said to be a \emph{thick
level}. It is worth emphasizing that the width of the interior
1-skeleton is a minimum taken up to isotopy rel $\partial M_i$.

\begin{figure}[h!]
\centering
\includegraphics{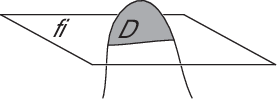}
\caption{An upper disc for $f_i$} \label{updisc}
\end{figure}

Consider a non-critical level $f_i$. Now consider a disc, $D$, with
the property that $\partial D$ consists of two arcs, one lying entirely on $f_i$ and the other running along an arc of 
interior 1-skeleton. Suppose also that the interior of $D$ is
disjoint from the interior 1-skeleton and that the disc emanates in
the upward (resp. downward) direction from $f_i$. Then $D$ is said
to be an \emph{upper} (resp. \emph{lower}) \emph{disc} for $f_i$.
Note that an upper or lower disc may intersect $f_i$ in its
interior. A simple example of an upper disc is shown in Figure \ref{updisc}.

If the interior 1-skeleton is in thin position with respect to $h$
then any pair of upper and lower discs for a non-critical level must
intersect at some point away from the interior 1-skeleton, for
otherwise we could reduce the overall width by isotoping across them
both. Another important thing to note is that a thick level always
has both an upper and a lower disc. For more on thin position, see
\cite{scharlthin}.

Our next goal is to show that the surfaces $K_i$ may be isotoped rel
boundary into interiorly-normal, interiorly-normal to one side or
almost interiorly-normal position. Our strategy will be to use a
similar inductive argument to that in \cite{stocking}. Much of this
argument is built on the following proposition. 

\begin{prop} \label{keyprop} Under the hypotheses of Theorem 3.1, each thick surface $K_i$ satisfies at least one fo the following:
\begin{enumerate}
  \item $K_i$ is isotopic to a surface in $M_i$ that is interiorly-normal,
  interiorly-normal to one side or almost interiorly-normal.
  \item $K_i$ is isotopic to a surface that is
  interiorly-normal or interiorly-normal to one side with essential
  tubes attached on one
  side. The interiorly-normal or interiorly-normal to one side surface
  is not normally parallel to $N_{i-1}$
  or $N_i$ on the side without the tubes attached. The side of the interiorly-normal
  or interiorly-normal to one side
surface which has the tubes attached is interiorly-normal.
\end{enumerate}
\end{prop}

During the proof of Proposition \ref{keyprop}, we will need to make use of the
following lemma, which closely resembles a fact about strongly
irreducible Heegaard splittings. See Lemma 6 of \cite{stocking}.

\begin{lemma} \label{keylemma}Suppose that $K_i \subseteq M_i$ is strongly irreducible
but compressible on both sides in $M_i$. Let $K_i'$ be the result of
compressing $K_i$ on one side in $M_i$. Then $K_i'$ is
incompressible on the other side.
\end{lemma}

\noindent \textbf{Proof.} Start by compressing $K_i$ on the side
away from $K_i'$ as much as possible, and call the resulting surface
$K_i''$. Let $X$ be the 3-manifold bounded by $K_i'$ and $K_i''$.
Suppose for a contradiction that there is a compression disc, $D$, for $K_i'$ on the side towards $K_i$.
Consider the intersection of $D$ with $K_i''$. Remove any trivial
innermost curves on $K_i''$ by an isotopy. If the resulting
intersection is now empty then $D$ is a compression disc for $K_i'$
in $X$. If the intersection is not empty then an innermost curve
must bound a compression disc for $K_i''$ in $X$, since $K_i''$ is
incompressible on the side away from $K_i$. In any case we may now
apply Lemma 3.5 of \cite{marionalpha} to obtain a compression disc, $D'$,
for $K_i'$ or $K_i''$ in $X$ which intersects $K_i$ in a single
essential curve. Now, since neither $K_i'$ nor $K_i''$ is parallel
to $K_i$, we have that $K_i$ is isotopic to a parallel copy of $K_i'$ or $K_i''$
with tubes attached. Furthermore, the third conclusion of Lemma 3.5 of \cite{marionalpha}  means we may take $D'$ to be disjoint from the tubes. The single simple closed curve of intersection of $D'$ with
$K_i$ cuts off a compression disc, $D''$, for $K_i$. This disc,
together with a meridian disc for one of the tubes, contradicts the
strong irreducibility of $K_i$. \hfill $\square$\\

\noindent \textbf{Proof of Proposition \ref{keyprop}.} There are three options:
\begin{enumerate}
  \item The interior 1-skeleton of $T$ has a thick level.
  \item The interior 1-skeleton of $T$ has no thick level, but there is at least one critical level.
  \item The interior 1-skeleton of $T$ has no critical level.
\end{enumerate}
We shall deal with each of these cases in turn:
\vskip 4pt
\noindent \textbf{Case 1: The interior 1-skeleton of $T$ has a thick
level.} In this case the first step is to find an interior fibre of
$h$ which intersects each face of the interior 2-skeleton in simple
closed curves and interiorly-normal arcs. Our method is very similar
to that in \cite{stocking} so we only give an outline here. The key
observation from that paper is that if an interior fibre of $h$
intersects the 2-skeleton in an arc which starts and ends at on same
edge, then there is an upper or lower disc for the level which lies
in the interior 2-skeleton. Now, just above the minimum at the
bottom of the thick region there must be a lower disc which lies in
the interior 2-skeleton (possibly after an isotopy of the interior
2-skeleton). Similarly, there must be an upper disc which lies in
the interior 2-skeleton just below the maximum at the top of the
thick region. Hence (after a small isotopy so that $h$ restricts to
a Morse function on each face of the interior 2-skeleton), at least
one of the following must be true:

\begin{enumerate}
  \item There is a level in the thick region with no upper or lower
  discs in the interior 2-skeleton.
  \item There is a level in the thick region with both upper and
  lower discs in the interior 2-skeleton.
  \item There is a level, $J$ say, in the thick region, which does not
  intersect the interior 2-skeleton in general position, with the property that a level just above
  it has an upper disc in the interior 2-skeleton and a level just below it
  has a lower disc in the interior 2-skeleton.
\end{enumerate}

  We may rule out option 2 straight away, since any pair of upper and
lower discs in the interior 2-skeleton must either fail to intersect
away from the interior 1-skeleton, violating thin position, or be
nested, again violating thin position. Option 3 can be ruled out by
noting that the pair of upper and lower discs for the levels just
above and below $J$ may be perturbed slightly to become a pair of
upper and lower discs for $J$ which fail to intersect away from the
interior 1-skeleton. The only option left is option 1, and so there
is a level, $L$ say, in the thick region which only intersects the
interior 2-skeleton in simple closed curves and interiorly-normal
arcs. This completes the first step.

For the second step, we aim to shown that $L$ may be isotoped rel 1-skeleton and
compressed on one side to obtain a new surface $L'$, where $L'$
intersects the interior 2-skeleton in interiorly-normal arcs and
$L'$ intersects each truncated tetrahedron in a collection of discs.

Consider a simple closed curve embedded on $L$ which bounds a disc
$D$, where $D$ lies entirely within the interior of a single ideal
tetrahedron $T$, the interior  of $D$ is disjoint from $L$, and
$\partial D$ does not bound a disc in $L \cap T$. We shall refer to
such a disc a \emph{local compression disc} for $L$. We start by
compressing and isotoping $L$ rel 1-skeleton so that it admits no
local compression discs, as follows. If $D$ is a genuine compression
disc for $L$ then compress along $D$. If not then $\partial D$
bounds a disc $D'$ in $L$ and $D \cup D'$ is a 2-sphere $S$.
Since $M$ is irreducible, $S$ bounds a 3-ball, and furthermore we
claim that $S$ does not intersect the interior 1-skeleton of the
ideal triangulation. To prove this suppose that $S$ intersects the
interior 1-skeleton. Let $N$ be a connected component of
intersection of $S$ with a truncated tetrahedron of the triangulation,
and suppose that $N$ intersects the interior 1-skeleton. Since $S$
is a 2-sphere, $N$ does not intersect the boundary faces of the
truncated tetrahedron. Also there is a component of $\partial N$ which
consists of normal arcs. Thus we may appeal to \cite{s3rec} for the
following facts about the length of each such component of $\partial
N$.

\begin{enumerate}
  \item If a component of $\partial N$ has odd length, then it
  has length 3.
  \item No component of $\partial N$ has length 6.
  \item If a component of $\partial N$ has length greater than 8 it
  crosses some edge of the interior 1-skeleton at least 3 times.
\end{enumerate}

If $\partial N$ has a component consisting of normal arcs of even
length, then it must have length at least 4, since a curve
consisting of just two normal arcs must consist of two arcs which
start and end on the same interior edge of the 1-skeleton. Now,
normal curves of length 3 bound triangles, normal curves of length 4
bound squares and normal curves of length 8 bound octagons.

We may rule out the possibility of a component of $\partial N$
intersecting an arc of interior 1-skeleton three times or more since
this would violate thin position, as in Claim 4.2 of \cite{s3rec}.
Hence each component of $\partial N$ which intersects the 1-skeleton
bounds a triangle, square or octagon. Note that this does not mean
that $N$ actually is a triangle, square or octagon, but it does mean
that we may replace $N$ with a collection of triangles, squares and
octagons without changing the part of the boundary of $N$ which
consists of normal arcs. Carry out this operation for all connected
components of intersection of $S$ with each truncated tetrahedron of the
triangulation which intersect the interior 1-skeleton, and throw
away the rest of $S$. Call the resulting surface $S'$. Then $S'$
must be a collection of 2-spheres, since it is homeomorphic to the
result of performing 2-surgery on a 2-sphere. But we assumed that
$M$ contains no 2-spheres consisting of just triangles, squares and
octagons. Hence $S$ does not intersect the interior 1-skeleton. Thus
we may isotope $L$ across the 3-ball which $S$ bounds without
changing the intersection of $L$ with the interior 1-skeleton.

Note that removing a local compression disc by compressing increases
the Euler characteristic of $L$ and removing one by an isotopy
reduces the number of intersections of $L$ with the interior
2-skeleton whilst not changing the Euler characteristic. Hence we
may isotope rel 1-skeleton and compress $L$ so that it admits no
local compression discs. Call the resulting surface $L'$. Since $L$
is strongly irreducible all the compressions were taken on the same
side. Suppose that a component of intersection of $L'$ with an ideal
tetrahedron is not a disc. Then we may find a compression disc for
that component. By an innermost curve argument, this compression
disc's interior may be assumed to be disjoint from $L'$ and hence it
is a local compression disc. This contradiction shows that $L'$
intersects each truncated tetrahedron in discs.

We claim that $L'$ intersects the interior 2-skeleton in
interiorly-normal arcs. Thus consider a simple closed curve of
intersection of $L'$ with a face of the interior 2-skeleton which is
innermost on that face. It bounds a disc $D$ in the face. By
pushing $D$, including its boundary, slightly into the neighbouring
tetrahedra in each direction, we cannot get a local compression
disc, and so $\partial D$ must bound a pair of discs in $L$, both of
whose interior is disjoint from the interior 2-skeleton. Together
these discs form a 2-sphere component of $L'$. But $L'$ has no
2-sphere components, a contradiction.

Hence we have achieved the second step, and $L$ may be compressed on
one side and isotoped rel 1-skeleton to obtain a surface $L'$, where
$L'$ intersects the interior 2-skeleton in interiorly-normal arcs
and $L'$ intersects each truncated tetrahedron in a collection of discs.

Our third step is to prove that $L'$ is either interiorly-normal, interiorly-normal to
one side or almost interiorly-normal (disc type). We shall use
arguments similar to those in \cite{s3rec}. First note that $L'$
intersects each ideal tetrahedron in a collection of discs, each of
which has boundary which intersects each interior edge at most
twice, for otherwise we would have a thinning pair of upper and
lower discs for $L$ (see claim 4.2 of \cite{s3rec}). Now, suppose
that $L'$ is not interiorly-normal or interiorly-normal to one side.
Then, since $L'$ is not interiorly-normal, $L'$ has a component of
intersection with an ideal tetrahedron which admits an edge
compression disc. Since $L'$ is not interiorly-normal to one side,
there must be a component of intersection of $L'$ with an ideal
tetrahedron which admits an edge compression disc on the other side.
But an edge compression disc for a component of intersection of $L'$
with an ideal tetrahedron is also an upper or lower disc for $L$,
and so any pair of edge compression discs for components of
intersection of $L'$ with ideal tetrahedra on opposite sides must
intersect away from the interior 1-skeleton. This cannot happen if
the two edge compression discs lie in different ideal tetrahedra,
and so all the edge compression discs must lie in the same ideal
tetrahedron. Call this ideal tetrahedron $\Gamma$. Suppose, for a
contradiction, that there are two edge compression discs for
components of intersection of $L'$ with $\Gamma$ on opposite sides
of $L'$ whose boundaries intersect different discs of intersection
of $L'$ with $\Gamma$. We will show that these edge compression
discs miss $L'$ in their interior. An innermost curve argument means
that we may suppose without loss of generality that these edge
compression discs have interiors which do not meet $L'$ in simple
closed curves. The edge compression discs cannot have interior
meeting $L'$ in arcs because this would establish the existence of
nested upper and lower discs for $L'$, violating thin position.
Hence the edge compression discs miss $L'$ in their interior.

But this means that the edge compression discs' interiors lie on
opposite sides of $L'$, and are hence disjoint. Hence the edge
compression discs intersect on $L'$, which contradicts the assertion
that their boundaries intersect different discs of intersection of
$L'$ with $\Gamma$. Thus any pair of edge compression discs for
components of intersection of $L'$ with $\Gamma$ on opposite sides
must be incident to the same disc of intersection of $L'$ with
$\Gamma$. This disc satisfies all the conditions to be almost
interiorly-normal and all the other discs of intersection of $L'$
with the ideal tetrahedra in the triangulation must be
interiorly-normal. Hence $L'$ must be almost interiorly-normal (disc
type) if it is not interiorly-normal or interiorly-normal to one
side. This completes the third step.

We now conclude the proof of Proposition \ref{keyprop} in this case. If no compressions were required when passing from $L$ to $L'$, then
option 1 of Proposition \ref{keyprop} holds. From now on assume that some
compressions were required. Now, if we start with $L'$ we may
recover $L$ by adding essential tubes dual to the compressions that
were carried out. Since $L$ is strongly irreducible, the tubes all
lie on the same side of $L'$ and, by Lemma \ref{keylemma}, this side is
incompressible. Call this side of $L'$ the $I$ side.

It is enough to show that we may isotope $L'$ rel boundary towards
the $I$ side so that it is interiorly-normal on that side. This is
guaranteed by the following lemma, whose proof is deferred to later,
and the following two observations. First, $L'$ has no 2-sphere
components. Second, $L'$ does not end up normally parallel to $N_i$
or $N_{i+1}$ on the side without the tubes attached because $L$ is a
thick level of $h$.

\begin{lemma} \label{pushlemma}Let $F$ be a properly embedded separating surface in $M_i$ that is
incompressible on one side, which we shall call the $I$ side.
Suppose $F$ satisfies all of the following properties:
\begin{enumerate}
  \item $F$ intersects each truncated tetrahedron in discs, apart from
  possibly one truncated tetrahedron which it intersects in discs and one
  annulus made out of two discs and a tube running parallel to an
  interior edge of the  1-skeleton which is attached on the $I$ side.
  \item $F$ intersects each face of the interior 2-skeleton in interiorly-normal
  arcs.
  \item $F$ intersects each face of the boundary 1-skeleton in
  normal arcs.
  \item $F$ admits at least one edge compression disc on the $I$
  side.
\end{enumerate}
Then $F$ may be isotoped rel boundary towards the $I$ side so that
each component is either interiorly normal on the $I$ side or a
2-sphere lying entirely within a truncated tetrahedron.
\end{lemma}
\noindent \textbf{Case 2: The interior 1-skeleton of $T$ has no thick level, but there is at least one critical level.}
In this case, since there is at least one critical level, there is
either at least one local maximum of the interior 1-skeleton or
there is at least one local minimum. Since there is no thick level,
no maximum can be at a greater height than a minimum. Hence all the
minima appear above all the maxima. We may suppose without loss of
generality that there are some minima, because otherwise there would
have to be some maxima and the proceeding proof will be similar.
Consider a level, $L$, just below the top of $M_i$. This level
consists of a surface parallel to $N_i$ which is interiorly-normal
in the downward facing direction, with tubes attached in the
downward direction. Since $L$ is above the top minimum, $L$ admits a
lower disc, $E$, the boundary of which consists of a sub-arc,
$\beta$, of an edge of the interior 1-skeleton and an arc, $\alpha$,
on $L$. The arguments we shall use are in large part identical to
those in Lemma 4 of \cite{stocking}, and we shall not repeat those
arguments in full here. Just as in \cite{stocking} define the
\emph{complexity} of $E$ to be $(a,b)$, where $a$ is the number of
intersections of the core graph of the tubes of $L$ with the
interior 2-skeleton, $b$ is the number of intersections of $E$ with
the interior 2-skeleton, and $(a,b)$ is ordered lexicographically.
We would like to isotope the tubes so that the complexity of $E$ is
$(0,0)$. Sadly this will not always be possible, but when the
arguments from \cite{stocking} break down we will be able to employ
a trick which saves the day. Consider a simple closed curve of
intersection of $E$ with the interior 2-skeleton which is innermost
on $E$. These may be removed via an isotopy as in Case 1 of the
argument in \cite{stocking}. Now consider an outermost arc of
intersection of $E$ with the interior 2-skeleton. The endpoints of
this arc both lie on $\alpha$. The reason that the arguments from
\cite{stocking} break down is that it is possible that the disc,
$D$, which this outermost arc cuts off of $E$ might not touch any
tubes at all. The good news is that, as we shall see, in this event
$D$ provides a recipe for isotoping $L$ in a useful way. For the
moment however, we shall proceed by using the arguments in
\cite{stocking} to remove outermost arcs of intersection of $E$ and
the interior 2-skeleton and thus reduce $L$'s complexity. Suppose
that at no stage do we see an outermost arc which cuts off a disc
which does not hit a tube. In this way we may reduce the complexity
of $E$ to $(0,0)$. There are a few points that should be made in
justification of using the arguments in \cite{stocking}, which
consists of consideration for several different cases.

\begin{figure}[h!]
\centering
\includegraphics{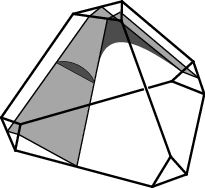}
\caption{A disc cut off by an outermost arc which hits no tubes}
\label{F:2}
\end{figure}

\begin{enumerate}
  \item In cases 3, 5, 7 and 8 in the proof of Lemma 4 of
  \cite{stocking} it is important that there cannot be a handle of
  $L$ contained in a tetrahedron. In the setting of \cite{stocking}
  one may appeal to Haken's Lemma, which does not apply here.
  Instead note that $L$ is a punctured 2-sphere, and so $L$ cannot have a handle contained
  inside a single truncated tetrahedron.
  \item Interiorly-normal (to the downward side) surfaces are different to the normal
  surfaces in \cite{stocking}, but this has no effect on our ability
  to transfer the arguments to this new setting, nor does the use of
  ideal triangulations.
  \item In case 8 of \cite{stocking} it is observed that the largest number of normal disks
  that can bound a connected component of a tetrahedron is four.
  The corresponding statement here is that there will only be a finite
  number of interiorly-normal or interiorly-normal to one side discs bounding a connected component
  of a truncated tetrahedron, although we do not know how big this number will be. This is enough for the rest of the argument in case 8 to be
  applied in this setting.
\end{enumerate}

Thus, provided that we never see an outermost arc which cuts off a
disc which does not hit a tube, we may reduce the complexity of $E$
to $(0,0)$. This means that we have performed an isotopy so that $E$
is contained in a single tetrahedron. We will now carry out an
isotopy within this tetrahedron so that $E$ runs over just one tube
which is parallel to an interior edge. This is achieved in exactly
the same way as in \cite{stocking}. If there are no other tubes then
option 1 holds. If there are some other tubes then they all lie on
the same side as the one that $E$ runs over. Compress the other
tubes so that the resulting surface is incompressible on the side
which contains $E$. Then option 2 holds by appealing to Lemma \ref{pushlemma}.

Now suppose that at some point in the above procedure we get an
outermost arc which cuts off a disc, $D$, which does not hit a tube.
Compress all of the tubes of $L$ and call the new surface $L'$. Then
$L'$ is incompressible on the same side as $D$, which we shall call
the $I$ side. Now, $D$ is a face compression disc for $L'$.

Note that since $L'$ is interiorly-normal on the $I$ side, it has no
edge compression discs on that side. Let us isotope $L'$ across $D$
and see what happens. Call the resulting surface $L''$. First note
that since $L'$ has no edge compression discs on the $I$ side, $L''$
intersects the all the truncated tetrahedra in discs which are
interiorly-normal on the $I$ side, apart from possibly those in the
tetrahedron opposite $D$. In this truncated tetrahedron, $T$, the effect
of isotoping $L'$ across $D$ is to band together two (possibly
non-distinct) interiorly-normal discs. If these discs are distinct
then $L''$ intersects $T$ in a collection of discs which are
interiorly-normal to the $I$ side and possibly one disc which admits
an edge compression disc on the $I$ side. If $L''$ is
interiorly-normal on the $I$ side then option 2 holds since $L''$ is
not normally parallel to $L'$. If $L''$ is not interiorly-normal to
the $I$ side then we may apply Lemma \ref{pushlemma}, and so option 2 holds in
this case as well.

Now suppose that the discs in $T$ that get banded together are not
distinct. Then $L''$ intersects $T$ in a collection of discs which
are interiorly-normal to the $I$ side, together with an annulus
which consists of two discs joined by a tube which lies on the
opposite side to the tubes which were compressed when passing from
$L$ to $L'$. Hence the local compression disc for $L''$
corresponding to this tube is not a genuine compression disc. By a
similar argument to that in case 1, this local compression disc must
cut off a 2-sphere which does not intersect the interior 1-skeleton. This contradicts the fact that $L'$ intersects the interior 2-skeleton in interiorly-normal arcs, meaning that this situation (where the discs in $T$ that get banded together are not
distinct) does not arise.

\vskip 4pt
\noindent \textbf{Case 3: The interior 1-skeleton of $T$ has no critical
level.} Consider a fibre of $h$ just below the top of $M_i$ and another one
just above the bottom. Both these levels have compression discs
which lie entirely within a truncated tetrahedron, but on opposite
sides. Hence (after a small isotopy so that $h$ restricts to a Morse
function on each face of the interior 2-skeleton) one of the
following must be true:

\begin{enumerate}
  \item There is a level which intersects the interior 2-skeleton in general position, with
  no compression disc whose boundary is contained in a single truncated tetrahedron.
  \item There is a level which intersects the interior 2-skeleton in general
  position, with compression discs
  on each side all of whose boundaries are contained in single truncated tetrahedra.
  \item There is a level which does not
  intersect the interior 2-skeleton in general position, with the property that a level just above
  it has on one side a compression disc whose boundary lies in a single ideal
  tetrahedron and a level just below it
  has on the other side a compression disc whose boundary lies in a single truncated tetrahedron.
\end{enumerate}

Suppose option 1 holds. Let $L$ be a level, which intersects the
interior 2-skeleton in general position, with no compression disc
whose boundary is contained in a single truncated tetrahedron. Hence any
local compression discs that $L$ admits must be non-essential, and
may therefore be removed with an isotopy rel 1-skeleton as in case
1. Let the result of removing all non-essential local compression
discs in this way be called $L'$. Suppose that $L'$ admits a
compression disc which lies entirely within a single truncated
tetrahedron. Then $L$ admits a compression disc whose boundary lies
entirely within a single truncated tetrahedron, a contradiction. Hence
$L'$ admits no local compression discs. Furthermore, $L'$ intersects
the interior 1-skeleton in interiorly-normal arcs, and since $L$ has
no edge compression discs, neither does $L'$. Hence $L'$ is
interiorly-normal and this proves Proposition \ref{keyprop} when option 1 holds.

We aim to rule out option 2. Suppose option 2 holds. If the
boundaries of the compression discs for the level, $L$ say, lie in
different truncated tetrahedra then they must be disjoint, contradicting
strong irreducibility. So they  lie in the same truncated tetrahedron,
$T$ say. Now, $L$ intersects the boundary of $T$ in a collection of
simple closed curves on $\partial T$. Consider a curve of $L \cap
\partial T$, which is innermost on $\partial T$ amongst those curves which
do not bound discs of $L \cap T$. This curve does not bound a disc
of $L \cap T$ but it does bound a disc whose interior lies in $T
\backslash L$. If this disc is a compression disc for $L$ then we
have contradicted the strong irreducibility of $L$. So it is
non-essential. Hence, as in case 1, we may apply an isotopy rel
1-skeleton to $L$ which reduces the number of components of
intersection of $L$ with the interior 2-skeleton to obtain a new
surface, $L'$ say. Note that $L'$ still has compression discs on
each side whose boundaries lie in $T$. We may now apply the same
argument to $L'$ and eventually we will contradict strong
irreducibility, ruling out option 2.

Suppose option 3 holds. Then there is a face of the interior
2-skeleton, $F$ say, which intersects the level, $L$ say, not in
arcs. If the truncated tetrahedra that this face bounds are distinct
then consider their union. Now remove from this a small open
neighbourhood of the faces of the interior 2-skeleton apart from
$F$. The result is topologically a 3-ball, so argue as in option 2
to contradict strong irreducibility.

\begin{figure}[h!]
\centering
\includegraphics{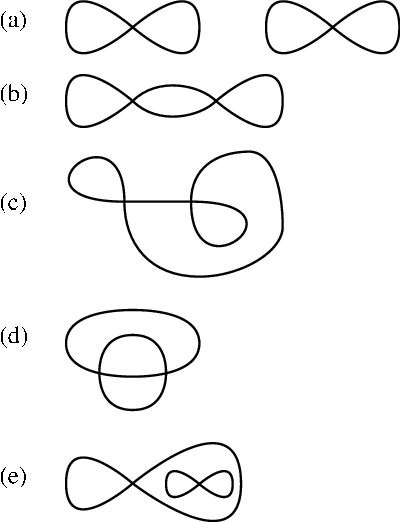}
\caption{The possible configurations of $L' \cap \partial T''$}
\label{last}
\end{figure}

Now suppose that the two tetrahedra which $F$ bounds are not
distinct. Let $T$ be the truncated tetrahedron which has two faces
identified to give the face $F$ of the interior 2-skeleton in $M$.
Consider $L \cap \partial T$, regarding $T$ as a 3-ball, and without
identified faces. Then $L \cap \partial T$ consists of a collection
of simple closed curves together with a graph which has just two
vertices, both with valance 4, as shown in Figure \ref{last}. If an
innermost simple closed curve of $L \cap
\partial T$ on $\partial T$ bounds a disc of $L \cap T$ then this disc, together with
a sub-disc of $\partial T$, bounds a 3-ball. Cut this 3-ball off
from $T$ and continue to cut off 3-balls in this manner as much as
possible. Call the result $T'$. Note that $T'$ is still
topologically a 3-ball and $L \cap
\partial T'$ still consists of simple closed curves and a graph as in Figure
\ref{last}. None of the simple closed curves bound discs of $L \cap
T'$, apart from possibly some curves which separate two components
of the 4-valent graph. Consider an innermost simple closed curve of
$L \cap
\partial T'$. It does not bound a disc of $L \cap T'$,
but it does bound a disc, $D$ say, whose interior lies in $T'
\backslash L$. Furthermore, $\partial D$ is disjoint from the two
compression discs of $L$ and so it must be non-essential for
otherwise we would contradict strong irreducibility. Hence $\partial
D$ bounds a sub-disc of $L$ which, together with $D$, form a
2-sphere. As in case 1, this 2-sphere must be disjoint from the
interior 1-skeleton, and so we may isotope across the 3-ball that
this 2-sphere bounds without affecting the intersection of $L$ with
the 1-skeleton. The isotopy has the same effect as performing
2-surgery along $D$ and throwing away the resulting 2-sphere
component. Let $S$ be the 2-sphere component that gets thrown away.
Now, we know that a level just above $L$ admits a compression disc,
$D'$ say, whose boundary lies entirely in $T'$ and runs along a band
in $L \cap T'$ which gets removed when we pass to a level just below
$L$. Without loss of generality we may suppose that $D'$ is disjoint
from $D$. Suppose that a 4-valent vertex of $L \cap
\partial T'$ appears on $S \cap
\partial T'$. Then the boundary of $D'$ lies on $S$. This
contradicts the fact that $D'$ is a compression disc for $L$. Hence
the isotopy rel 1-skeleton does not affect the 4-valent graph part
of $L \cap
\partial T'$. After
having performed the isotopy rel 1-skeleton we may have a new
component of $L \cap
\partial T'$ which bounds a disc of $L \cap T'$. If so, then use
this disc to cut off another 3-ball from $T'$. Continue to reduce
the number of simple closed curves of $L \cap
\partial T'$ in this manner until there are none, apart from
possibly some simple closed curves separating two components of the
4-valent graph. Let the resulting sub-manifold of $T$ be called
$T''$ and the result of isotoping $L$ be called $L'$. Then $L' \cap
\partial T''$ consists of a graph, as shown in Figure
\ref{last}, and in case (a) or (e) possibly some simple closed
curves separating the components of the graph.

We wish to show that in configurations (a) and (e), $L' \cap
\partial T''$ may be arranged to contain no simple closed curve components. Suppose
that $L' \cap
\partial T''$ consists of a graph as in configuration (a) or (e) and
possibly some simple closed curves which separate the two components
of the graph. Suppose that $L' \cap \partial T''$ does have a simple
closed curve component, $c$ say. Suppose that $c$ does not bound a
disc of $L' \cap T''$. Then $L' \cap T''$ admits a local compression
disc, $D$, for $L'$ in $T''$. Because $L$ is as described in option
3, $D$ must be non-essential. Hence there is a sub-disc, $D'$, of
$L'$, whose boundary agrees with that of $D$ and so that $D \cup D'$
forms a 2-sphere bounding a 3-ball disjoint from the interior
1-skeleton. Hence we may isotope rel 1-skeleton across the 3-ball.
As before, this isotopy does not affect the 4-valent graph part of
$L' \cap
\partial T''$. Perform these isotopies rel 1-skeleton as much as
possible and call the resulting surface $L''$. Then all simple
closed curve components of $L'' \cap \partial T''$ must bound discs
of $L'' \cap T''$.

We seek to show that $L'' \cap \partial T''$ has no simple closed
curve components. Suppose the contrary. Then the simple closed curve
components of $L'' \cap \partial T''$ bound a collection of parallel
discs in $L'' \cap T''$. Now, these discs separate separate $T''$
into a number of components. Two of these components have a 4-valent
graph component of $L'' \cap \partial T''$ on their boundary. Call
them $T_1$ and $T_2$. Let $L''_+$ be a level just above $L''$. Note
that for $i = 1,2$,  $L''_+ \cap \partial T_i$ consists of either a
single simple closed curve or two simple closed curves. If $L''_+
\cap \partial T_i$ consists of a single simple closed curve then
$L''_+ \cap T_i$ consists of a disc since $L''$ is planar. If $L''_+
\cap \partial T_i$ consists of two simple closed curves then, since
$L''$ is planar, $L''_+ \cap T_i$ consists of an annulus or a pair
of discs. Hence either $L''_+ \cap T_1$ or $L''_+ \cap T_2$ consists
of an annulus, because otherwise $L''_+$ would have no compression
discs whose boundary lies in $T''$. Without loss of generality
suppose that $L''_+ \cap
\partial T_1$ is an annulus. Now, let $L''_-$ be the level just below
$L''$. When we pass from $L''_+$ to $L''_-$ we isotope across a face
compression disc which intersects the co-core of the annulus in
$T_1$ just once. The effect which this has in $T_2$ is to add a band
to $L''_+ \cap \partial T_2$. Hence $L''_- \cap \partial T_2$
consists of an annulus. This annulus admits a face compression disc
which intersects the co-core of $L''_- \cap \partial T_2$ just once.
Hence $L''$ admits a pair of compression discs which intersect just
once, namely at the identified 4-valent vertex of $L'' \cap
\partial T''$. This is a contradiction. Hence $L'' \cap \partial T''$ has no simple closed
curve components.

Suppose that $L'' \cap \partial T''$ is a graph as in configuration
(a). Remember that as we pass from a level just above $L''$ to a
level just below $L''$ we are isotoping the entire surface in the
same direction. Hence there is a level, $L'''$ say, which is either
just above or just below $L''$ and which intersects $\partial T''$
in two simple closed curves. They cannot bound discs because $L$ is
as described in option 3. Hence they bound an annulus because $L$ is
planar. Now, there must be a face compression disc for $L'''$ in
$T''$ which hits the co-core of the annulus. We isotope across this
face compression disc when passing to the level the other side of
$L''$. But that means that the graph of $L'' \cap
\partial T''$ should be connected, a contradiction.

Suppose that $L'' \cap \partial T''$ is a graph as in configuration
(b). Again remember that as we pass from a level just above $L''$ to
a level just below $L''$ we are isotoping the entire surface in the
same direction. Hence there is a level either just above or just
below $L''$ which intersects $\partial T''$ in a single simple
closed curve. Since $L$ is planar, it must be a disc. But this
contradicts the fact that $L$ is as described in option 3.

Now suppose that $L'' \cap \partial T''$ is a graph as in
configuration (c) or (d). Then a level, $L'''$ say, either just
above or just below $L''$ intersects $\partial T''$ in two simple
closed curves. As in configuration (a), these curves must bound an
annulus. Also, the annulus must admit a face compression disc which
hits its co-core. Now, $L'''$ intersects the interior 2-skeleton in
interiorly-normal arcs and simple closed curves. If $L'''$
intersects $\partial T''$ in a simple closed curve on the interior
2-skeleton, disjoint from the 1-skeleton, then that means that $L'''
\cap \partial T''$ consists of two such curves, one each on the two
faces which get identified when forming $M$. But that means that
$L'''$ has a torus component, which it does not. Hence $L'''$
intersects $\partial T''$ in interiorly-normal arcs. Furthermore,
since $L$ has no edge compression discs, $L'''$ intersects $T''$ in
an annulus made of two interiorly-normal discs joined by a face
parallel tube. Moreover, $L'''$ has no compression discs contained
in any tetrahedra other than $T$. Hence, after an isotopy rel
1-skeleton to remove any non-essential local compression discs,
$L'''$ may be ambient isotoped to be almost interiorly-normal.

Suppose that $L'' \cap \partial T''$ is a graph as in configuration
(e). Then a level, $L'''$ say, either just above or just below $L''$
intersects $\partial T''$ in three simple closed curves. Since $L$
is planar, $L''' \cap T''$ must consist of either a three times
punctured 2-sphere or an annulus and a disc. In the later case we
may argue as in configurations (c) and (d). So suppose that $L'''
\cap T''$ consists of a three times punctured 2-sphere. Suppose that
$L''' \cap T''$ admits a compression disc in $T''$ which is not a
genuine compression disc for $L'''$ in $M$. Then remove this local
compression disc with an isotopy rel 1-skeleton as in case 1. The
resulting surface intersects $T''$ in an annulus. If its co-core is
not essential then we have found a surface as described in option 1.
So suppose that the co-core is essential. Then we may argue as in
configurations (c) and (d). Hence we may suppose that all the
compression discs for $L''' \cap T''$ in $T''$ are genuine
compression discs for $L'''$. They all lie on the same side of
$L'''$ because otherwise we could argue as in option 2. Now, as we
pass from $L'''$ to the level the other side of $L''$, we isotope
across a face compression disc. This has the effect of boundary
compressing $L'''$ at the same time as adding a band, when
considering $L'''$ as a properly embedded surface in $T''$. Hence
the surface, $L''''$ say, the other side of $L''$ to $L'''$ admits a
compression disc on the same side as those of $L'''$ in $T''$ whose
boundary lies in a single truncated tetrahedron. But we know that
$L''''$ admits a compression disc on the other side whose boundary
lies in a single truncated tetrahedron. This means we are in option 2, a
contradiction. \hfill $\square$
\vskip 4pt
We deferred the proof of Lemma \ref{pushlemma} during the proof of Proposition \ref{keyprop}. We rectify this now.

\begin{namedtheoremb}  Let $F$ be a properly embedded separating surface in $M_i$ that is
incompressible on one side, which we shall call the $I$ side.
Suppose $F$ satisfies all of the following properties:
\begin{enumerate}
  \item $F$ intersects each truncated tetrahedron in discs, apart from
  possibly one truncated tetrahedron which it intersects in discs and one
  annulus made out of two discs and a tube running parallel to an
  interior edge of the  1-skeleton which is attached on the $I$ side.
  \item $F$ intersects each face of the interior 2-skeleton in interiorly-normal
  arcs.
  \item $F$ intersects each face of the boundary 1-skeleton in
  normal arcs.
  \item $F$ admits at least one edge compression disc on the $I$
  side.
\end{enumerate}
Then $F$ may be isotoped rel boundary towards the $I$ side so that
each component is either interiorly normal on the $I$ side or a
2-sphere lying entirely within a truncated tetrahedron.
\end{namedtheoremb}

\noindent \textbf{Proof.} Our strategy will be to isotope across
edge compression discs on the $I$ side and remove local compression
discs on the $I$ side with an isotopy. If $F$ intersects an ideal
tetrahedron in an annulus then push the tube so that it surrounds
the edge  it was parallel to. Continue by isotoping across edge
compression discs on the $I$ side. If at any stage $F$ admits a
local compression disc, then we claim that this disc lies on the $I$
side. Suppose, on the contrary, that at some point we first have a
component of intersection, $C$ say, of $F$ with a tetrahedron, $T$,
which is an annulus (or possibly a surface of even lower Euler
characteristic) which compresses on the non-$I$ side. To return to
the previous step we must isotope towards the non-$I$ side. The
effect of this on $C$, considered as a sub-manifold of the 3-ball
$T$ is either to band together two points on its boundary, or to
boundary compress towards the non-$I$ side. Neither of these
operations can have the effect of returning it to being a disc, and
so no local compression discs appearing within a tetrahedron appear
on the non-$I$ side. Hence every local compression disc that appears
is on the $I$ side, and these may be removed with an isotopy rel
interior 1-skeleton as in case 1 of the proof of Proposition 3.4.
Isotoping across an edge compression disc reduces interior edge
degree and removing local compression discs decreases the number of
components of intersection of $F$ with the interior 2-skeleton
without affecting the interior
edge degree. Hence this process terminates. Each component of the resulting surface satisfies all the requirements to be interiorly-normal to the $I$ side provided it is not a 2-sphere lying entirely within a single truncated tetrahedron. \hfill $\square$ \\

We apply Proposition \ref{keyprop} to prove Theorem \ref{nktheorem} as follows. If option 1
holds, then we  stop. If option 2 holds then  we may remove
tubes from a surface isotopic to $K_i$ and obtain a surface which is
interiorly-normal to the resulting incompressible side and which is
not parallel to $N_{i-1}$ or $N_i$ on the side without the tubes
attached. Call this surface $N'$. Now cut $M_i$ along $N'$ and throw
away everything to the side to which the tubes are not attached. Now apply Proposition \ref{keyprop} again, and continue to do so repeatedly. Note that if any
stage of this iteration yields a surface which is boundary parallel
on the side with the tubes attached then the next application of
Proposition \ref{keyprop} will be via case 3 which in turn yields option 1. We
claim that if this process is repeated eventually option 1 will
hold, and this is what we shall now prove. Suppose, on the contrary,
that when we apply Proposition \ref{keyprop} in this manner option 2 holds
indefinitely. Let $N^{(1)} = N'$ and for $i > 1$ let $N^{(i)}$ be
the sequence of surfaces yielded by Proposition \ref{keyprop}. To get a
contradiction we would like to say that there can only be finitely
many non-parallel, disjoint interiorly-normal or interiorly-normal
to one side surfaces in $M$. This, however, is not true, for
consider an infinite sequence of non-parallel boundary parallel
interiorly-normal annuli. This example illustrates the extra
information that we have about the surfaces $N^{(i)}$, namely that
since the surfaces $N^{(i)}$ are related by isotopies rel boundary,
their boundaries are parallel as normal curves on $\partial M$. Let
$b$ be the boundary edge degree of the surfaces $N^{(i)}$.

For each $i$, consider the intersection of $N^{(i)}$ with each edge
of the boundary 1-skeleton of $M$. This will constitute a sequence
of points along that edge. For any two of these points, either they
are joined by an arc of intersection of $N^{(i)}$ with the interior
2-skeleton or they are not. There are only finitely many choices as
to which pairs of points on the same boundary edge are joined in
this manner. Hence we may find a subsequence $N^{(i_k)}$ of
$N^{(i)}$ where each surface in the subsequence carries the same
information as to which points of intersection with each edge of the
boundary 1-skeleton are joined to each other by an arc on the
interior 1-skeleton. Let $N = \cup_{k=1}^r  N^{(i_k)}$, where $r$ is
arbitrarily large.

The following is inspired in part by the proof of Lemma 13.2 of
\cite{hemp}. With this in mind, observe that a truncated tetrahedron of
$M$ is cut into pieces by $N$. We define a \emph{bad piece} as one
which contains a \emph{bad point}, which we now describe. Consider a
face of the 2-skeleton, which may be either an interior face or a
boundary face. Note that $N$ intersects the face in finitely many
collections of parallel copies of normal or interiorly-normal arcs.
Place a bad point in the interior of each component of the face
which is not bounded by two parallel normal or interiorly-normal
arcs. Let $n$ be the number of bad points in a given ideal
tetrahedron, $T$. Consider the normal or interiorly-normal arcs on
the boundary of $T$ and remove all the interiorly-normal arcs which
run from a boundary edge to the same boundary edge. The number of
bad points that this removes cannot exceed $b$. Hence the remaining
number of bad points, $n'$ say, is at least $n - b$. But there are at
most 56 remaining bad points, namely 4 on each boundary face and 10
on each interior face. Hence $56 \geq n' \geq n-b$ and so $n \leq 56
+ b$. Hence the number of bad pieces is at most $56 + b$.

Note that $\beta_1 (M) < \infty$ and so $M$ admits at most finitely
many compact properly embedded incompressible surfaces whose union
is non-separating. That means that $N$ cuts $M$ into arbitrarily
many pieces, for large enough $r$. Hence there must be a component
$C$ of $M \backslash N$ that contains no bad point, for large
enough $r$. Since $M$ is orientable, the closure of $C$ must be a
product bundle, and $C$ bounds two interiorly-normal (possibly to
one side) surfaces which are parallel as interiorly-normal (possibly
to one side) surfaces. This is a contradiction and Theorem \ref{nktheorem} is proved. \hfill $\square$

\vskip 4pt

It is worth remembering that Theorem \ref{nktheorem} assumes that $M$ contains no
2-spheres made out of just triangles, squares and octagons. We shall
see why this assumption is justified in the  Section 4.

\section{The Algorithm}

We will now turn our attention to the main theorem of this paper,
namely that there exists an algorithm to search for bridge punctured
2-spheres of given Euler characteristic for $M$, the exterior of a hyperbolic knot in $S^3$. An overview of our
algorithm to do this proceeds as follows:
\vskip 4pt
  \noindent \textbf{Step 1:} Construct a suitable ideal triangulation for the
  knot exterior, together with some extra information about how the
  boundary edges look relative to the natural foliation of
  the boundary by meridians.
  \vskip 4pt
  \noindent \textbf{Step 2:} Amongst the infinitely many types of
  interiorly-normal, interiorly-normal to one side and almost
  interiorly-normal discs, construct a finite subset of these types out of which we may build
  the thin and thick surfaces of $\mathcal{B}'$, the generalised bridge
  surface as in Theorem 2.1.
  \vskip 4pt
  \noindent \textbf{Step 3:} Write down the system of matching equations for
  the disc types found in Step 2. Algorithmically solve these
  equations to find a finite collection of fundamental solutions.
  \vskip 4pt
  \noindent \textbf{Step 4:} Use the fundamental solutions found in Step 3 to
  specify a finite list of candidates for $\mathcal{B}'$.
  \vskip 4pt
  \noindent \textbf{Step 5:} Algorithmically amalgamate each candidate for $\mathcal{B}'$, to obtain
  a finite list of candidate bridge punctured 2-spheres.
  \vskip 4pt
  \noindent \textbf{Step 6:} For each candidate bridge punctured 2-sphere, inspect it to
  see if it is indeed a bridge decomposition for the exterior of $K$ corresponding to a bridge punctured 2-sphere
  with the correct Euler characteristic.
\vskip 4pt
The remainder of this paper will be devoted to further developing
this overview.

In \cite{tn1alg} Marc Lackenby introduced the notion of a partially
flat angled ideal triangulation and used them to exhibit an
algorithm to determine the tunnel number of a hyperbolic knot in
$S^3$. We will now describe these triangulations in more detail. The following definitions apply to either a genuine ideal triangulation for a non-compact 3-manifold or to the sort that we consider for compact 3-manifolds, built out of truncated tetrahedra. 

An \emph{angle structure} for an ideally triangulated 3-manifold is
an assignment of interior angles in the range $(0,\pi)$ to each
interior edge of each tetrahedron in the triangulation so that the
angles associated to opposite edges are equal, the angles around
each ideal vertex sum to $\pi$ and the sum of the angles around each
edge is $2\pi$. If we allow some of the angles to be either $0$ or
$\pi$ then there may be some \emph{flat ideal tetrahedra}, which
consist of an angle of $\pi$ assigned to one pair of opposite edges
and angles of $0$ assigned to the other two pairs of opposite edges.
Pairs of faces which cobound an edge with interior angle $\pi$ are
said to be \emph{coherent}.

\begin{figure}[h!]
\centering
\includegraphics{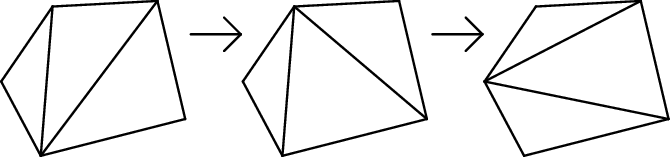}
\caption{Elementary moves on an ideally triangulated ideal pentagon}
\label{movonpoly}
\end{figure}

A \emph{layered polygon} is a collection of flat ideal tetrahedra
glued together in a certain way. They are defined in \cite{tn1alg}
as arising from a sequence of elementary moves applied to an ideal
polygon with ideal triangulation as follows. Start with an ideal
polygon with ideal triangulation, and suppose that we apply to it a
sequence of elementary moves to change the triangulation in such a
way that every edge in the interior of the ideal polygon is affected
by an elementary move. An example of such a sequence of elementary moves is
shown in Figure \ref{movonpoly}.

\begin{figure}[h!]
\centering
\includegraphics{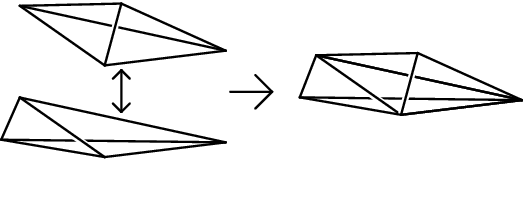}
\caption{Building a layered polygon} \label{laypoly}
\end{figure}

A layered polygon arises from such a sequence of elementary moves as follows.
For each move take a flat ideal polygon and place two coherent faces
either side of the edge that is removed in that move. This gives a
new ideal triangulation for the ideal polygon. Continue by placing
more flat ideal tetrahedra underneath for each elementary move, as
shown
in the Figure \ref{laypoly}.

A layered polygon built in this way has a special type of edge on
its boundary; namely those on the boundary of the original ideal
polygon. These edges are called the \emph{vertical boundary} of the
layered polygon.

A \emph{partially flat angled ideal triangulation} is an ideal
triangulation, with a real number in the range $[0, \pi]$ assigned
to each edge of each ideal tetredhron, satisfying the following
conditions:

\begin{enumerate}
  \item The angles at each ideal vertex of each ideal tetrahedron sum to
  $\pi$.
  \item The angles around each edge sum to $2 \pi$.
  \item If the angles of a tetrahedron are not all strictly positive, then the tetrahedron is
flat.
  \item The union of the flat tetrahedra is a collection of
layered polygons, possibly with some edges in their vertical
boundary identified.
\end{enumerate}

The following theorem appears as Theorem 2.1 in \cite{tn1alg}, and
is absolutely key in that paper.

\begin{theorem}[Theorem 2.1 of \cite{tn1alg}]\label{21lack} Any finite-volume hyperbolic 3-manifold $M$ with
non-empty boundary has a partially flat angled ideal triangulation.
Moreover, there is an algorithm that constructs one, starting with
any triangulation of
$M$.
\end{theorem}

In \cite{tn1alg}, Lackenby searches for surfaces which do not
intersect the boundary of the 3-manifold. This is not the case when
searching for bridge surfaces, and we will need to put more control
on the behavior of the boundary of the partially flat angled ideal
triangulations with which we wish to work. Our goal is to
algorithmically find a partially flat angled ideal triangulation for
the knot exterior together with an upper bound on the number of
intersections with the boundary 1-skeleton that are required for a
meridian.

Let $M$ be a knot exterior. Let $\partial M = S^1_M \times S^1_L$ be
the product structure on $\partial M$ by meridians and longitudes.
As stated in Section 2, we know from \cite{straighttorus} that the
boundary 1-skeleton of an ideal triangulation may be isotoped so
that the following conditions hold:

\begin{enumerate}
  \item All the boundary edges are transverse to the foliation of $\partial
  M$ by meridional circles.
  \item All the vertices of the triangulation have different
  meridional coordinates.
\end{enumerate}

When these two conditions hold we shall say that the boundary
1-skeleton is in \emph{standard position}. Note that there may be
many different ways of placing the boundary 1-skeleton in standard
position. We shall refer to a meridional leaf of the foliation of
$\partial M$ by meridional circles which intersects a vertex of the
boundary 1-skeleton as a \emph{singular meridian}. Other meridians
are \emph{non-singular}.

The \emph{edge degree} of a non-singular meridian is the number of
times that it intersects the boundary 1-skeleton. The
\emph{meridional edge degree} of a standard position of a
triangulation of $\partial M$ is the maximum edge degree of all the
non-singular meridians. The \emph{minimal meridional edge degree} of
a triangulation of $\partial M$ is the minimal meridional edge
degree of the boundary 1-skeleton taken over all isotopies of the
triangulation into standard position.

The following theorem represents the first step in our algorithm to
search for bridge punctured 2-spheres of given Euler characteristic.

\begin{theorem} \label{specialideal}Let $M$ be the exterior of a hyperbolic knot,
$K$, in $S^3$. Starting with a diagram of $K$, there exists an
algorithm to construct a partially flat angled ideal triangulation
for $M$ together with an upper bound on the minimal meridional edge
degree of the ideal triangulation.
\end{theorem}

\noindent \textbf{Proof of Theorem \ref{specialideal}.} We know from \cite{tn1alg} that there
exists a partially flat angled ideal triangulation for $M$. To find
one algorithmically starting from any ideal triangulation, Lackenby
argues as follows. There is an algorithm to test whether an ideal
triangulation admits a partially flat angle structure, and this is
simply a linear programming question. We also know from Theorem
1.2.5 of \cite{mat} that any two ideal triangulations are connected
by a sequence of 2-3 and 3-2 moves. So, to find a partially flat angled ideal triangulation of $M$
starting with any ideal triangulation, $T$, for $M$ we test it to
see if it admits a partially flat angle structure. If it does not
then apply all possible 2-3 and 3-2 moves to the ideal triangulation
and test all the resulting ideal triangulations in the same manner.
Eventually we must find an ideal triangulation, $T'$, which does
have a partially flat angle structure and this is where the
algorithm stops. This partially flat angled ideal triangulation is
the one for which we wish to place an upper bound on its minimal
meridional edge degree. To achieve this we need to keep track of the
meridians of the knot throughout.

Start with a knot diagram for $K$ and remove any nugatory crossings
so that resulting diagram is reduced. It is a theorem of Menasco
(see \cite{tn1}) that any reduced diagram canonically induces an
ideal polyhedral decomposition of the knot exterior with just two
ideal polyhedra. Construct this ideal polyhedral decomposition and
mark a meridian of the knot exterior on the ideal boundary.
Subdivide the decomposition so as to obtain an ideal triangulation,
$T$, for the knot exterior, and keep track of the meridian on the
boundary.

\begin{figure}[h!]
\centering
\includegraphics{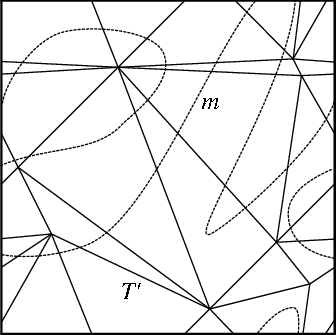}
\caption{$\partial M$} \label{sccthereon}
\end{figure}

Now apply 2-3 and 3-2 moves to $T$ to obtain $T'$, the ideal
triangulation of the knot exterior which we know has a partially
flat angle structure. Still keep track of the meridian on the
boundary of the ideal triangulation. We wish to calculate an upper
bound on the minimal meridional edge degree of $T'$, and we will use
the meridian, $m$, to help us. The boundary of $T'$ is simply a
triangulated torus, and $m$ is a simple closed curve thereon, as
shown in the Figure \ref{sccthereon}.

\begin{figure}[h!]
\centering
\includegraphics{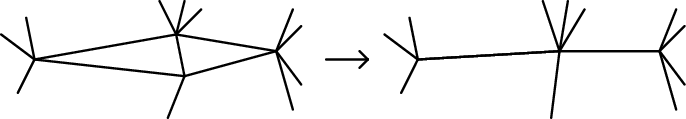}
\caption{A contraction of an edge} \label{contrac}
\end{figure}

We now need to take a closer look at Bojan Mohar's
\cite{straighttorus} proof that any simple triangulation of a
torus may be ambient isotoped so that all its edges are geodesic
line segments. A \emph{contraction} of an edge of a triangulation is
the move shown in Figure \ref{contrac}.

The strategy that Mohar uses in \cite{straighttorus} is to apply as
many contractions to the given triangulation as possible. When no
more contractions are possible he shows that the resulting
triangulation must be homeomorphic to that shown in Figure
\ref{standardtri} below. This triangulation is clearly isotopic to
one with geodesic edges. To recreate the original triangulation we
carry out the inverse operation to contraction in a small
neighbourhood of the edges which have been changed.

\begin{figure}[h!]
\centering
\includegraphics{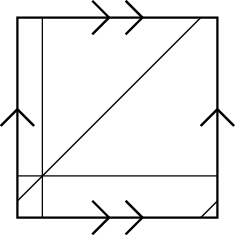}
\caption{A triangulation of the torus with just three edges}
\label{standardtri}
\end{figure}

\begin{figure}[h!]
\centering
\includegraphics{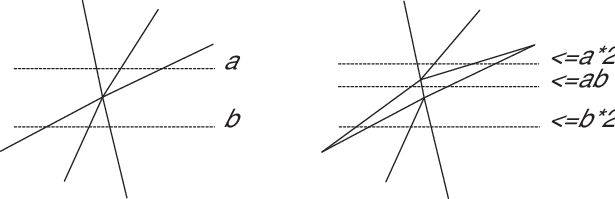}
\caption{Undoing the contractions} \label{undocont}
\end{figure}

With Mohar's ideas in mind, apply all possible contractions to $T'$,
and keep track of $m$. Now, $m$ might not intersect the new
triangulation in normal arcs. Rectify this by isotoping it to remove
any non-normal arcs. The resulting triangulation, $T''$ intersects
$m$ in $n$ points, say. Note that the minimal meridional edge degree
of $T''$ is at most $n$. Now carry out the inverse operation to
contraction in a small neighbourhood of the edges which were changed
in passing from $T'$ to $T''$. Each time we apply an inverse
contraction to $T''$, we create a new vertex in the triangulation
and three new edges, as shown in Figure \ref{undocont}.

With each inverse contraction we create one more non-singular
meridian. This has edge degree at most double the total of the
meridional edge degree of the non-singular meridians to either side
before the inverse contraction. Furthermore, any other non-singular
meridian's edge degree is increased at most by a factor of 2. Hence
the effect of an inverse contraction on the minimal meridional edge
degree is to increase it by at most a factor of 4. Hence the minimal
meridional edge degree of $T'$ is certainly at most $4^{t}n$ where
$t$ is the number of edges of $T'$. \hfill $\square$ \\

We shall now return to the task proving Theorem 1.1. The following theorem is an empowered version of Theorem \ref{nktheorem}.

\begin{theorem}\label{newthm} Let $K$ be a knot is $S^3$. Let $M$ be the exterior of $K$ and let $M$ have an ideal triangulation which contains no 2-spheres consisting of just triangles, squares and octagons. Let $F$ be bridge punctured $2$-sphere for $M$. Let $\mathcal{B'}$ be a generalized bridge surface for $M$ such that:
\begin{enumerate}
  \item Each thin surface of $\mathcal{B}'$ is incompressible.
  \item Each thick surface of $\mathcal{B'}$, $K_i$, is strongly
irreducible in $M_i$, the submanifold of $M$ obtained by cutting along the adjacent thin surfaces $N_i$ and $N_{i-1}$.
  \item No thick surface, $K_i$, cobounds a product $(\textrm{Surface}) \times
  I$ with an adjacent thin surface $N_i$ or $N_{i-1}$.
  \item  $\mathcal{B'}$ amalgamates to give $F$.
\end{enumerate}
Then there is a computable constant $b$, which may be calculated
from the Euler characteristic of $F$ and the the minimal meridional edge degree of the triangulation
of $\partial M$, such that $\mathcal{B'}$ may be ambient isotoped so that:
\begin{enumerate}
  \item The thin surfaces of $\mathcal{B'}$ are interiorly normal.
  \item The thick surfaces of $\mathcal{B'}$ are interiorly-normal, interiorly-normal to one
side or almost interiorly-normal.
\item $\mathcal{B'}$ intersects the boundary 1-skeleton of $T$ at most $b$ times.
  \end{enumerate}
\end{theorem}

\noindent \textbf{Proof of Theorem \ref{newthm}.} Let $u$ be the minimal meridional edge degree of the triangulation of $\partial M$. We know that the
surfaces of $\mathcal{B}'$ each have meridional boundary. Hence
these surfaces are isotopic to ones whose boundary components
intersect the boundary 1-skeleton at most $u$ times and which
intersect the boundary 2-skeleton in normal arcs. Now, we can
calculate an upper bound on the number of surfaces $N_i$ and $K_i$
of $\mathcal{B}'$ from the Euler characteristic of $F$. We also know that each of these surfaces has
no more boundary components than $F$, which in turn has $2-\chi(F)$
boundary components. Hence we know an upper bound on the total
number of boundary components of $\mathcal{B}'$. Hence we have an upper bound on the boundary edge degree of $\mathcal{B}'$ at this stage. Now apply Theorem 3.1 to isotope $\mathcal{B}'$ so that it satisfies the first two requirements of the conclusion of Theorem 4.3. Note that our bound on the boundary edge degree of $\mathcal{B}'$ is preserved under this isotopy because the boundary of $\mathcal{B}'$ is only changed with a boundary height adjusting isotopy. Thus it is this bound that we should take as $b$. \hfill $\square$\\

\noindent \textbf{Proof of Theorem \ref{maintheorem}.} Let $K$ be a hyperbolic knot in $S^3$ and let $M$ be the exterior of $K$. The algorithm begins by finding the ideal
triangulation, $T$, of Theorem \ref{specialideal}, together with an upper bound,
$u$, on its minimal meridional edge degree. Note that  one of the
hypotheses for Theorem \ref{newthm} is that $M$ contains no embedded
2-spheres consisting of just triangles, squares and octagons. This
fact is due to the existence of a partially flat angle structure on
$T$, our ideal triangulation for $M$. See Theorem 2.2 of
\cite{tn1alg}. Thus, by Theorems \ref{untelprop} and  4.3, for any bridge punctured 2-sphere $F$ of $M$ there is a generalized bridge decomposition
$\mathcal{B}'$ that satisfies the conclusions of Theorem \ref{newthm}.

The upper bound on boundary edge degree $b$ gives us
an upper bound on the number of times one of the interiorly-normal,
interiorly-normal to one side or almost interiorly-normal discs
which make up $\mathcal{B}'$ intersect the boundary 1-skeleton.
Now, it is clear that there are only finitely many
interiorly-normal, interiorly-normal to one side or almost
interiorly-normal disc types which intersect the boundary edges at
most a given number of times. Furthermore they may be found
algorithmically since there will be only finitely many paths
consisting of normal and interiorly-normal arcs on the boundary of a
given truncated tetrahedron to inspect.

Thus we know that every thick or thin surface of $\mathcal{B}'$ is made up of a patched together collection of finitely
many different disc types (and possibly one face parallel tube).
Furthermore these discs may be found algorithmically. Suppose that
each truncated tetrahedron admits $d$ disc types from this collection,
and that there are $t$ truncated tetrahedra in the ideal triangulation
of $M$ found in Theorem \ref{specialideal}. Then each of the interiorly-normal,
interiorly-normal to one side or almost interiorly-normal surfaces
in $M$ specify a vector in $V = (\mathbb{N} \cup 0)^{dt}$, where
each coordinate represents the number of each different disc type in
the surface. The vector representing a surface, $S$, will be denoted
by $f(S)$. This representation of a surface by a vector is in the
same spirit as classical normal surface theory.

Note that on each interior face of the ideal triangulation the discs
of an interiorly-normal, interiorly-normal to one side or almost
interiorly-normal surface must patch together. Now, each disc type
in two neighbouring truncated tetrahedra gives rise to a number of
interiorly-normal arcs on the interior face which joins them. The
number of interiorly-normal arcs of each type on this face arising
from the discs in each of the two tetrahedra must be the same in
order for them to patch together to form a surface. Thus for a
vector in $V$ to represent a surface, the vector must satisfy a
system of linear equations, known as the \emph{matching equations}.
Note that the matching equations are specified by an ideal
triangulation as well as a finite collection of disc types.

We will return to the subject of matching equations shortly, but
first we will consider how a solution to the matching equations may
give rise to a surface. The essential fact is that a solution, $v
\in V$, to the matching equations may correspond to no surface, or
it may correspond to one or more surfaces, but in any case it is
possible to algorithmically find all surfaces, $S$, such that $f(S)
= v$.

Returning to the matching equations, consider two vectors, $v_1, v_2
\in V$ which satisfy these equations. Then $v_1 + v_2$ is also a
solution to the matching equations. A solution to the matching
equations which cannot be written in the form $v_1 + v_2$ for
non-zero solutions $v_1$ and $v_2$ is called a \emph{fundamental
solution} to the matching equations. The following theorem is key, for which the reader is referred to  \cite{hemion} for a proof.

\begin{theorem} There exists an algorithm to calculate all the
fundamental solutions to the system of matching equations.
\end{theorem}

So far we have not made great use of the fact that the ideal
triangulation for $M$ has a partially flat angle structure. A
surface $F$ in general position with respect to the triangulation
inherits a \emph{combinatorial area} (sometimes referred to as just \emph{area}) from a partially flat angle
structure as follows. Let $T$ be a truncated tetrahedron and consider a
connected component $D$ of $F \cap T$. Note that $\partial D$ hits
the edges of $\partial T$ transversely. The \emph{area} of $D$ is defined
to be the sum of the exterior angles of these edges, counted
with multiplicity, minus $2\pi$ times the Euler characteristic of
$D$, where the \emph{exterior angle} of an edge is taken to be $\pi$ minus
the interior angle. The interior angle at a boundary edge of an
truncated tetrahedron is deemed to be $\frac{\pi}{2}$. The combinatorial area of $F$
is the sum of the areas of all the components of intersection of $F$
with each truncated tetrahedron. A simple calculation implies that the
combinatorial area of $F$ is $-2\pi\chi(F)$. Hence if $F$ is a 2-sphere then it
has negative area. But a quick check tells us that triangles,
squares and octagons have non-negative area. Hence there are no
2-spheres made of just triangles, squares and octagons in an ideal
triangulated 3-manifold with a partially flat angle structure. One
can in fact go further. The following theorem is due to Lackenby and
appears as Theorem 2.2 in \cite{tn1alg}.

\begin{theorem}[Theorem 2.2 of  \cite{tn1alg}] \label{paratori}Let $T$ be a partially flat angled ideal
triangulation of $M$. Then any connected 2-normal surface in $T$
with non-negative Euler characteristic is normally parallel to a
boundary
component.
\end{theorem}

A 2-normal surface is one which consists of triangles, squares and
octagons. We shall refer to an embedded surface made up of
interiorly-normal, interiorly-normal to one side and almost
interiorly-normal discs as \emph{interiorly 2-normal}. Note that an
interiorly 2-normal surface which does not intersect the boundary is
necessarily 2-normal.

Suppose that $v_1, \ldots, v_m \in V$ are the fundamental solutions
to the matching equations. Let $S$ be the union of all the surfaces
of $\mathcal{B}'$, where each component of $S$ is interiorly-normal,
interiorly-normal to one side or almost interiorly-normal. Then
$f(S) = \sum_{i=1}^{n}a_{i}{v_i}$ for non-negative integers $a_i$.
In order to find a finite list of candidates for $f(S)$ we need to
bound each $a_i$. This is achieved in the proof of the following
theorem, which we defer to the very end of this paper.

\begin{theorem} \label{findthem}Let $T$ be a partially flat angled ideal
triangulation of $M$. Then, for any positive integers $n$ and $b$,
$T$ contains only finitely many properly embedded surfaces $F$ such
that:
\begin{enumerate}
  \item $- \chi(F) \leq n$
  \item The boundary edge degree of $F$ is at most $b$
  \item Each component is either interiorly-normal,
  interiorly-normal to one side or almost interiorly-normal
\end{enumerate}
Furthermore, there is an algorithm to construct all of these
surfaces.
\end{theorem}

Theorem \ref{findthem} means that we may algorithmically find a finite
list of candidates for $\mathcal{B}'$, a generalized bridge
surface obtained by applying Theorems 2.1 and \ref{newthm} to a bridge punctured 2-sphere of given Euler characteristic for a given hyperbolic knot in $S^3$. Note that the Euler characteristic of $\mathcal{B}'$ may be estimated from that of the bridge punctured 2-sphere.  There is an algorithm to determine
whether $\mathcal{B}'$ is a generalised bridge surface. This is
achieved in a similar way to Section 5, Step 3 of \cite{tn1alg}. We
have therefore reduced the task at hand to one of inspecting each
candidate for $\mathcal{B}'$ to see if it amalgamates to give a
bridge decomposition for $M$ with the required genus. The task of
algorithmically amalgamating is achieved in essentially the same way
as in \cite{tn1alg}, and we shall not repeat the details here.
Finally, by Theorem 4.1.13 of \cite{mat} there is an algorithm to
test whether a separating punctured 2-sphere for $M$ which
intersects $\partial M$ in meridians is a bridge punctured 2-sphere.
This completes the proof of Theorem \ref{maintheorem}, assuming Theorem \ref{findthem}. \hfill $\square$\\

We complete this paper by proving Theorem \ref{findthem}.\\

\noindent \textbf{Proof of Theorem \ref{findthem}.} The first step of the algorithm is to
write down the matching equations corresponding to the ideal
triangulation and the interiorly 2-normal disc types which intersect
the boundary at most $b$ times. Now algorithmically solve these
equations to obtain a finite list of vectors $v_1, \ldots, v_m \in
V$ which are the fundamental solutions to these matching equations.
For a surface $F$ which satisfies the conditions in the Theorem, we
know that $$f(F) = \sum_{i=1}^{m}a_{i}{v_i}$$ for non-negative
integers $a_i$. Recall that we may algorithmically find all surfaces
which correspond to a given vector. We have therefore reduced the
task at hand to one of bounding the magnitude of the integers $a_i$.
To achieve this, start by noting that we can associate to each
solution vector $v$ of the matching equations its boundary edge
degree $b(v)$ which is simply the sum of the number of times the
disc types of $v$ intersect the boundary. Hence,
$$\sum_{i=1}^{m}a_{i}b({v_i}) \leq b .$$ Now, for those $i$ such
that $b(v_i) > 0$, we have $a_i \leq \frac{b}{b(v_i)} =: b_i$.
Without loss of generality suppose that the vectors with $b(v_i) >
0$ are precisely those with $1 \leq i \leq p \leq m$ for some $p$.

\begin{figure}[h!]
\centering
\includegraphics{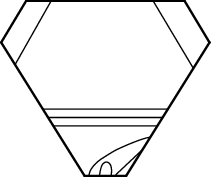}
\caption{Boundary parallel edges in an interior face are parallel}
\label{paraedges}
\end{figure}

Let $$V_1 = \sum_{i=1}^{p}a_{i}{v_i}$$ and $$V_2 =
\sum_{i={p+1}}^{m}a_{i}{v_i}.$$ Then $V_2$ is a solution to the
matching equations for which every non-zero coordinate corresponds
to a disc which does not intersect any boundary arc of the
1-skeleton. Hence $V_2$ represents a closed embedded 2-normal
surface. We now wish to show that $V_1$ represents a properly
embedded interiorly 2-normal surface. To achieve this, we start by
noting that on a given interior face of the 2-skeleton of $M$ all
the arcs of intersection of $F$ with that face with endpoints on the
same pair of interior arcs of the 1-skeleton must be parallel on
that face. Call these edges the \emph{boundary parallel edges}. They
are illustrated in Figure \ref{paraedges}.

\begin{figure}[h!]
\centering
\includegraphics{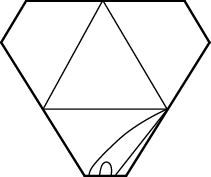}
\caption{The intersection of $C(F)$ with an interior face of the two
skeleton} \label{paraedgessq}
\end{figure}

Now form a 2-complex, which we shall denote $C(F)$, by homotoping
the boundary parallel edges together, as well as well as all points
of intersection of $F$ with an interior arc of the 1-skeleton, as
shown in Figure \ref{paraedgessq}.

To form a surface which is represented by $V_1$, simply delete from
$C(F)$ one of each disc corresponding to $V_2$ and tease the
remaining 2-complex apart along the boundary parallel arcs in its
1-skeleton. Let $f(F_1) = V_1$ and for $i = p+1,\ldots,m$ let
$f(S_i) = v_i$. Then $$f(F) = f(F_1) +
\sum_{i={p+1}}^{m}a_{i}f(S_i)$$ and so $$\chi(F) = \chi(F_1) +
\sum_{i={p+1}}^{m}a_{i}\chi(S_i).$$ Hence $$\chi(F_1) +
\sum_{i={p+1}}^{m}a_{i}\chi(S_i) \geq -n.$$ Hence $$|\chi(F_1)| +
\sum_{i={p+1}}^{m}a_{i}\chi(S_i) \geq -n.$$ Therefore
$$\sum_{i={p+1}}^{m}a_{i}\chi(S_i) \geq -n - |\chi(F_1)|.$$ Now, we know from Theorem \ref{paratori} that for $i = p+1, \ldots,
m$ either $S_i$ has negative Euler characteristic or it is a
boundary parallel torus. Let $v_m$ represent the boundary parallel
torus. Hence for $i = p+1, \ldots, m-1$ we have $$a_i \leq n +
|\chi(F_1)|.$$ But we have already reduced $V_1$ to one of finitely
many possibilities, and so $F_1$ is one of only finitely many
candidates. Hence we can find un upper bound for $|\chi(F_1)|$, and
consequently an upper bound for $a_i$ when $i = p+1, \ldots, m-1$.
This leaves the issue of boundary parallel tori. For this we need to
make use of the following claim.
\vskip 4pt

\begin{claimx} Let $f(F) = v + v_m$ where $v_m$ represents a
boundary parallel normal torus. Then there exists an interiorly
2-normal surface $F'$ such that $f(F') = v$ and $F$ is obtained from
$F'$ by performing a switch along the collection of simple closed
curves of intersection of $F'$ with a boundary parallel torus.
\end{claimx}
\noindent \textbf{Proof of Claim.} As before form the 2-complex $C(F)$.
Now, rather than deleting all the triangles corresponding to $v_m$,
isotope them a little towards the boundary of $M$ so that we have a
boundary parallel torus, $R$. Denote the remainder of $C(F)$ as
$C(F)'$.

By an innermost curve argument, the triangles of $R$ intersect
$C(F)'$ in arcs. Furthermore these arcs must start and end on the
intersection of a normal arc of $R$ and an arc joining a boundary
edge and a neighbouring interior edge. Tease apart the 2-complex
$C(F)'$ to obtain $F'$.

\begin{figure}[h!]
\centering
\includegraphics{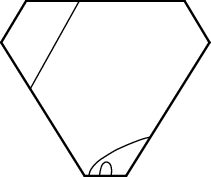}
\caption{The arc types of the discs of $v$} \label{3typ}
\end{figure}

Now, since $F$ contains a triangle of every type in each ideal
tetrahedron, the discs corresponding to non-zero coordinates of $v$
cannot intersect the interior 2-skeleton in arcs running from a
boundary edge to another boundary edge, or arcs running from a
boundary edge to the opposite interior edge. This means that the
discs corresponding to non-zero coordinates of $v$ must intersect
the interior 2-skeleton in arcs which either run from one interior
edge to another, or from a boundary edge to the same boundary edge,
or from a boundary edge to a neighbouring interior edge. These arc
types are illustrated in Figure \ref{3typ}. Furthermore, the discs
corresponding to non-zero coordinates of $v$ cannot intersect two
different boundary triangles of the 2-skeleton of an ideal
tetrahedron.

\begin{figure}[h!]
\centering
\includegraphics{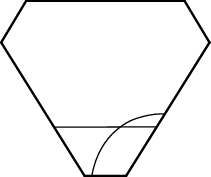}
\caption{Arcs of $R$ and $F'$ on an interior face} \label{int}
\end{figure}

Consider an interior face of the 2-skeleton. If components of
intersection of $R$ and $F'$ with this face intersect then it must
be as shown in Figure \ref{int}.

\begin{figure}[h!]
\centering
\includegraphics{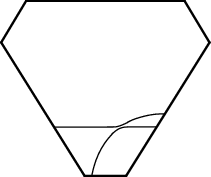}
\caption{A regular switch} \label{intas}
\end{figure}

\begin{figure}[h!]
\centering
\includegraphics{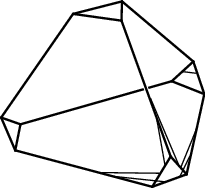}
\caption{Six arc types of $F'$ which might hit $P$} \label{sixarcs}
\end{figure}

If we wish to perform a switch to resolve this point of
intersection, then it must be as shown in Figure \ref{intas}.
We shall call such a switch a \emph{regular switch}. Consider the
point of intersection before performing a regular switch. This is
the endpoint for an arc of intersection of $F'$ and $R$ in an ideal
tetrahedron on each side. We wish to show that $R$ and $F'$ may be
ambient isotoped rel 2-skeleton so that switches may be performed
along all curves of intersection in such a way that the switches are
regular on every face of the interior 2-skeleton. Let $T$ be an
truncated tetrahedron of $M$. Consider a connected component of $R \cap
T$, which we shall call $P$. There are six types of arc of $F' \cap
\partial T$ on the boundary of this truncated tetrahedron
which might intersect $P$, as shown in Figure \ref{sixarcs}.

\begin{figure}[h!]
\centering
\includegraphics{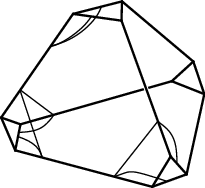}
\caption{Three configurations for $F'$} \label{3discs}
\end{figure}

Three of these arc types turn to the left when moving away from the
boundary of $M$ and three turn to the right. For regular switches to
agree along an arc of intersection of $F'$ and $R$ in $T$ they must
run between a left turning arc and a right turning arc. Suppose $P
\cap F' \neq \emptyset$. Now, it cannot be the case that $F' \cap
\partial T$ has only left turning arcs or only right turning arcs
emanating from a given boundary face. Hence there is at least one
left turning arc and at least one right turning arc emanating from
the boundary face which $P$ is parallel to. Consider the first left
turning arc to the right of a right turning arc and the right
turning arc to the left of this. This pair of arcs must appear in
one of the three configurations shown in Figure \ref{3discs}.

Isotope $F'$ rel 2-skeleton so that there is an arc of intersection
of $F'$ and $P$ running between these two arcs, and perform a
regular switch along this arc. Now look at the remaining left
turning and right turning arcs and apply a similar argument.
Eventually we will have resolved all arcs of intersection of $F'$
and $P$ in $T$ in a way which is regular on $\partial T$. Carry out
this procedure in each tetrahedron to complete the proof of the
claim. \hfill{$\square$}

We will say that $F$ is obtained from $F'$ by \emph{adding a
boundary parallel torus}. We will slightly abuse notation and write
$F = F' + R$. If it were the case that all of the regular switches
took place around essential arcs on $R$ and furthermore these
switches were oriented in the same direction around $R$, then $F'$
and $F$ would be ambient isotopic. This, however, need not be the
case. For this reason we will need to take a closer look at what
happens if boundary parallel tori are repeatedly added to a surface.

Remember that the task at hand is to bound $a_m$. With this in mind,
note that by repeatedly applying the previous claim we arrive at a
surface $F''$ whose $v_m$ coordinate is zero and for which $F = F''
+ a_m.R$. Consider the surfaces $F'' + k.R$ where $k$ is a
non-negative integer. These surfaces are made by performing regular
switches on the the disjoint union of $F''$ and a collection of $k$
parallel copies of $R$. Let $S_1$ and $S_0$ be tori, $S_1$ just
above the top (furthest away from the boundary) copy of $R$ and
$S_0$ just below the bottom one. Then $S_1$ and $S_0$ form the
boundary of an $I$-bundle, $S \times I$, with $S \times \{i\} = S_i$
for $i = 1,0$.

\begin{figure}[h!]
\centering
\includegraphics{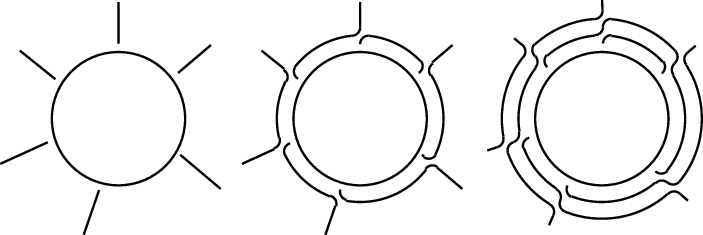}
\caption{Adding boundary parallel tori which switch along essential
cures} \label{swtch}
\end{figure}

Now, $F'' + k.R$ intersects $S \times I$ in a collection of properly
embedded surfaces. The boundary of these surfaces is a collection of
curves on $S_1$ and $S_0$. Let $c$ be the number of these curves on
each of $S_1$ and $S_0$. The collection of curves on each of $S_1$
and $S_0$ is the same when they are projected onto $S$. Note that
the essential curves must all be parallel. We'll turn our attention
to the non-essential curves shortly, but for now suppose there are
none. A schematic of this situation is shown in Figure \ref{swtch}.
Observe that $F'' + k.R$ intersects $S \times I$ in a collection of
annuli. The boundary of one such annulus consists of two essential
curves on $S_1 \cup S_0$. Note that if two curves on $S_1$ (resp.
$S_0$) are joined by an annulus of intersection of $F'' + k.R$ with
$S \times I$ then they still are for $F'' + k'.R$ when $k'>k$.
Suppose that two curves on $S_1$ (resp. $S_0$) eventually get joined
by an annulus of intersection of $F'' + k.R$ with $S \times I$ for
large enough $k$. Let the number of essential curves between these
two curves be $d$ (that is, if one projects the annulus onto $S_1$
(resp. $S_0$) then $d$ is the number of essential curves which lie
in the interior of this projection). Clearly $d<c$. Then these two
curves get joined by an annulus of intersection of $F'' + k.R$ with
$S \times I$ whenever $k > \frac{d}{2}$. (Note $d$ is always even.)
This is because an annulus such as this only intersects at most the top
(resp. bottom) $\frac{d}{2} + 1$ copies of $R$. Let $d'$ denote the
maximum value of $d$ taken over all pairs of essential curves on
$S_1$ (resp. $S_0$) which eventually get joined by an annulus of
intersection of $F'' + k.R$ with $S \times I$ for large enough $k$.
Let $k > 2(\frac{d'}{2}+1) + 1=d'+3$. This holds when $k > c + 3$.
When this is the case, there is a copy of $R$ in $F'' + k.R$ such
that, after switching, every part of it is contained in an annulus
joining opposite sides of $S \times I$.

Now consider the essential curves that do not eventually get joined
to another curve on the same side of $S \times I$. There are $e\leq
c$, say, of these curves on each of $S_1$ and $S_0$. Switches along
these curves are all oriented the same way. When $k > c+3$, adding
extra tori has the effect of making each of these curves on $S_1$
join to the next one around on $S_0$. This is not necessarily an
isotopy, but when $k
>c+3+e$ any $F'' + k.R$ must be related to $F'' + (k-e).R$ by an
isotopy. Hence for the case when there are no non-essential curves
on $S_1$ or $S_0$ we should restrict $a_m$ to at most $2c+3$.

Now suppose that there are some non-essential curves on $S_1$ and
$S_0$. Let $x$ be the total number of these curves. Consider the
disjoint union of $F''$ and the $k$ parallel copies of $R$ in $S
\times I$. Let us switch along the non-essential curves of
intersection only, and throw away the part of $F''$ in $S \times I$
whose boundary components are essential on $S_1$ and $S_0$. The
result is a properly embedded surface in $S \times I$, which we
shall call $F'' +_n k.R$. Now, a non-essential curve on $S_1$ or
$S_0$ can only be connected to the top or bottom $x$ copies of $R$
in $F'' +_n k.R$. Hence if $k > 2x+1$ then $F'' +_n k.R$ contains a
properly embedded torus, $R'$. We may therefore take a new product
neighbourhood of $R'$ which intersects $F''$ only in essential
curves on $R'$. Hence we may now argue as in the case where there
are no non-essential curves to conclude that if
$k>4c+4\geq(2x+1)+(2c+3)$ then $F''+k.R$ is isotopic to
$F''+(k-e).R$. Hence we should restrict $a_m$ to $4c+4$.

In order to estimate $c$, note that each essential curve along which
a switch is carried out must hit an arc of $F'$ which runs from a
boundary edge to a neighbouring interior edge of the 1-skeleton.
Hence $c$ is certainly at most $b$. The proof of Theorem \ref{findthem} is
completed by restricting $a_m$ to at most $4b+4$. \hfill$\square$
\vskip4pt

\bibliographystyle{amsplain}
\bibliography{bridge}

\providecommand{\bysame}{\leavevmode\hbox to3em{\hrulefill}\thinspace}
\providecommand{\MR}{\relax\ifhmode\unskip\space\fi MR }
\providecommand{\MRhref}[2]{%
  \href{http://www.ams.org/mathscinet-getitem?mr=#1}{#2}
}
\providecommand{\href}[2]{#2}
\begin{thebibliography}{10}

\bibitem{bachderbsedg}
Davic Bachman, Ryan Derby-Talbot, and Eric Sedgwick, \emph{Almost normal
  surfaces with boundary}, 2012.

\bibitem{bach}
David Bachman, \emph{Heegaard {S}plitting with {B}oundary and {A}lmost {N}ormal
  {S}urfaces}, Topology and its Applications \textbf{116} (2001), 153--184.

\bibitem{marionalpha}
Marion~Moore Campisi, \emph{$\alpha$-sloped generalized heegaard splittings of
  3-manifolds}, 2011.

\bibitem{ep}
David Epstein and Robert Penner, \emph{Euclidean decompositions of noncompact
  hyperbolic manifolds}, J. Diff. Geom. \textbf{27} (1988), 67--80.

\bibitem{gabaithin}
David Gabai, \emph{Foliations and the topology of 3-manifolds {I}{I}{I}}, J.
  Diff. Geometry \textbf{26} (1987), 479--536.

\bibitem{hsh}
Chuichiro Hayashi and Koya Shimokawa, \emph{Thin position of a pair
  (3-manifold, 1-submanifold)}, Pacific J. Math. \textbf{197} (2001), 301--324.

\bibitem{hemion}
Geoffrey Hemion, \emph{The classification of knots and 3-dimensional spaces},
  Oxford Science Publications, 1992.

\bibitem{hemp}
John Hempel, \emph{3-manifolds}, Annals of Mathematics Studies, Princeton Univ.
  Press, 1976.

\bibitem{normtori}
William~H. Jaco and J.~Hyam Rubinstein, \emph{0-efficient triangulations of
  three-manifolds}, J. Diff. Geometry \textbf{65} (2003), 61--168.

\bibitem{heeggenuslack}
Marc Lackenby, \emph{The {H}eegaard genus of amalgamated 3-manifolds}, Geom.
  Dedicata \textbf{109} (2004), 139--145.

\bibitem{tn1}
\bysame, \emph{Classification of alternating knots with tunnel number one},
  Comm. Anal. Geom. \textbf{13} (2005), 151--186.

\bibitem{tn1alg}
\bysame, \emph{An algorithm to determine the {H}eegaard genus of simple
  3-manifolds with non-empty boundary}, 2007.

\bibitem{mat}
Sergei Matveev, \emph{Algorithmic {T}opology and {C}lassification of
  3-{M}anifolds}, Algorithms and Computation in Mathematics, vol.~9, Springer,
  2003.

\bibitem{straighttorus}
Bojan Mohar, \emph{Straight-line representations of maps on the torus and other
  fiat surfaces}, Discrete Mathematics \textbf{155} (1996), 173--181.

\bibitem{scharlthin}
Martin Scharlemann, \emph{Handbook of {K}not theory}, ch.~Thin position in the
  theory of classical knots, pp.~429--459, Elsevier, 2005.

\bibitem{untel}
Martin Scharlemann and Abigail Thompson, \emph{Thin position for 3-manifolds},
  AMS Contemporary Math. \textbf{164} (1994), 231--238.

\bibitem{stocking}
Michelle Stocking, \emph{Almost normal surfaces in 3-manifolds}, Trans. Amer.
  Math. Soc. \textbf{352} (2000), 171--207.

\bibitem{s3rec}
Abigail Thompson, \emph{Thin {P}osition and the {R}ecognition {P}roblem for
  ${S}^3$}, Mathematical Research Letters \textbf{1} (1994), 613--630.

\bibitem{robinbridge}
Robin Wilson, \emph{Meridional {A}lmost {N}ormal {S}urfaces in {K}not
  {C}omplements}, 2007.

\end{thebibliography}

\end{document}